\documentclass[10pt,a4paper]{article}

\usepackage{amssymb}
\usepackage{amsmath}
\usepackage{tensor}
\usepackage{mathrsfs}
\newtheorem{Th}{Theorem}
\newtheorem{Prop}{Proposition}[section]
\newtheorem{Co}{Corollary}

\newtheorem{Lma}{Lemma}[section]

\newtheorem{Rm}{Remark}

\newcommand{\be}{\begin{equation}}
\newcommand{\ee}{\end{equation}}
\newcommand{\bes}{\begin{equation*}}
\newcommand{\ees}{\end{equation*}}
\newcommand{\R}{\mathbb{R}}

\newcommand{\C}{\mathbb{C}}

\newcommand\res{\mathop{\hbox{\vrule height 7pt width .5pt depth 0pt
\vrule height .5pt width 6pt depth 0pt}}\nolimits}

\def\theequation{\thesection.\arabic{equation}}
\def\theTh{\Roman{section}.\arabic{Th}}
\def\theProp{\Roman{section}.\arabic{Prop}}
\def\theCo{\Roman{section}.\arabic{Co}}

\def\theRm{\Roman{section}.\arabic{Rm}}
\newcommand{\reset}{\setcounter{equation}{0}\setcounter{Th}{0}\setcounter{Prop}{0}\setcounter{Co}{0}\setcounter{Lma}{0}\setcounter{Rm}{0}}

\def\al{\alpha}
\def\la{\lambda}
\def\eps{\varepsilon}

\def\Om{D_1(0)\setminus\{0\}}

\def\bn{\vec{n}}

\def\bex{\bAe_1}
\def\bey{\bAe_2}

\def\bei{\bAe_i}
\def\bej{\bAe_j}
\def\bek{\bAe_k}

\def\px{\partial_{x_1}}
\def\py{\partial_{x_2}}
\def\pj{\partial_{x^j}}

\def\pz{\partial_{z}}
\def\pzb{\partial_{\bar{z}}}
\def\bAe{\vec{e}}
\def\bH{\vec{H}}

\def\bA{\vec{A}}
\def\bF{\vec{F}}
\def\bU{\vec{U}}
\def\bJ{\vec{J}}

\def\bL{\vec{L}}

\def\bR{\vec{R}}
\def\bE{\vec{E}}
\def\bX{\vec{X}}
\def\bV{\vec{V}}
\def\bPe{\vec{P}}
\def\bw{\vec{w}}
\def\bu{\vec{u}}

\def\bB{\vec{B}}

\def\bp{\vec{\Phi}}
\def\bP{\vec{\Phi}}
\def\bxi{\vec{\xi}}

\def\bT{\vec{T}}
\def\bQ{\vec{Q}}

\def\bul{\bullet}
\def\di{D_1(0)}
\newcommand{\D}[3]{\tensor*[^{#1}]{D}{}_{#2}^{#3}}
\newcommand{\gammi}[3]{\tensor*[^{#1}]{\Gamma}{}_{#2}^{#3}}
\newcommand{\nabli}[2]{\tensor*[^{#1}]{\nabla}{}_{#2}}

\def\res{\mathop{\hbox{\vrule height 7pt width .5pt 
depth 0pt\vrule height .5pt width 6pt depth 0pt}}\nolimits}

\begin{document}


\reset

\title{Ends of Immersed Minimal and Willmore Surfaces\\
 in Asymptotically Flat Spaces}
\author{Yann Bernard\footnote{School of Mathematical Sciences, Monash University,
3800 Victoria, Australia.}\:\:,\:Tristan Rivi\`ere\footnote{Department of Mathematics, ETH Zentrum,
8093 Z\"urich, Switzerland.}}
\maketitle
\noindent
$\textbf{Abstract:}$ {\it We study ends of an oriented, immersed, non-compact, complete Willmore surfaces, which are critical points of the integral of the square of the mean curvature, in asymptotically flat spaces of any dimension; assuming the surface has $L^2$-bounded second fundamental form and satisfies a weak power growth on the area. We give the precise asymptotic behavior of an end of such a surface. This asymptotic information is very much dependent on the way the ambient metric decays to the Euclidean one. Our results apply in particular to minimal surfaces.}

\section{Introduction}

\subsection{Setting and Main Results}

Let $m\ge3$ be an integer, and let $(M,h_M)$ be a smooth and complete Riemannian manifold of dimension $m$. We will suppose that $(M,h_M)$ is asymptotically flat, i.e. that there exists a compact set $Z\subset M$ such that $M\setminus Z$ consists of finitely many ends, namely $M\setminus Z=\bigcup_{k=1}^{N}E_k$. Each end $E_k$ is diffeomorphic to $\mathbb{R}^m\setminus B^m_{r_k}(0)$, where $B^m_{r_k}(0)\subset\mathbb{R}^m$ is the ball of radius $r_k>0$ centered at the origin. Let $f_k:E_k\rightarrow\mathbb{R}^m\setminus B^m_{r_k}(0)$ be this diffeomorphism. Let $p$ denote the asymptotically flat coordinate induced by $f_k$. We require that the pull-back metric satisfy (for each $k$)
\bes\label{defAF}
h_{\al\beta}(p)\;:=\;\Big(\big(f_k^{-1}\big)^*h_M\Big)_{\al\beta}(p)\;=\;\delta_{\al\beta}+b_{\al\beta}(p)\:,
\ees
with\footnote{Throughout this paper, we will use the following standard notation. We write $f(X)=\text{O}_N(|X|^s)$ to indicate that $f^{(j)}(X)=\text{O}(|X|^{s-j})$ for all integers $j\in[0,N]$.}
\be\label{condAF}
b_{\al\beta}(p)\;=\;\text{O}_2(|p|^{-\tau})\qquad\text{for some}\:\:0<\tau\leq1 \quad\text{and for}\:\:|p|\,\gg\,1\:.
\ee
We will henceforth assume that $M\setminus Z$ has only one such end, diffeomorphic, say, to $\R^m\setminus B^m_1(0)$. \\[-1ex]

In the literature, the asymptotic behavior of the remainder $b_{\al\beta}(p)$ is dictated by the applications which one has in mind. Oftentimes, the metric is chosen to decay to the Schwarzschild metric (a so-called ``strongly asymptotically flat" condition) \cite{Car1, Car2,  HY, LMS, Met}. This essentially amounts to choosing $\tau=1$, along with some radial-only dependency condition on $b_{\al\beta}$. This Schwarzschild-type hypothesis ultimately follows from the proof of the positive mass theorem by Schoen and Yau \cite{SY}.  When questions of stability are involved \cite{Car1, Car2, HY}, one must require that the decay hold to four derivatives. In this paper, we will only be concerned with obtaining information on the second fundamental form, which is why we only require the asymptotic behavior of $b_{\al\beta}$ to hold up to second-order derivatives. In \cite{Hua}, where the existence of a foliation by constant mean curvature spheres is shown, the author requires the asymptotic decay to satisfy a so-called Regge-Teitelboim condition, namely (\ref{condAF}) with $1\ge\tau>1/2$. The present work is concerned, in parts, with finding results that can hold for the smallest possible value $\tau$. \\

We study a certain class of complete non-compact surfaces in $(\R^m,h)$, namely Willmore surfaces, which will be made precise below. Let us point out that minimal surfaces are Willmore surfaces, so all of our results apply in particular to complete non-compact minimal surfaces in asymptotically flat space. \\ 

Let $S$ be a connected, oriented, non-compact, complete, two-dimensional surface immersed in $(\R^m,h)$. We let $\vec{A}_S^h$ denote the second fundamental form of $S$ (this is a normal vector mapping into $\R^m$, hence the arrow notation).  We assume that
\be\label{teague}
\int_S|\bA^h_S|_h^2\,d\mu_h\;<\;\infty\:,
\ee
where $\mu_h$ denote the induced measure on $S$. One can also understand $S$ as a complete immersed surface into $\R^m$ equipped with the Euclidean metric. Naturally, the corresponding fundamental form $\bA^{h_0}_S$ differs from $\bA^h_S$. One might wonder whether (\ref{teague}) holds with the Euclidean metric $h_0$ in place of $h$. This is in general false, and an additional hypothesis is needed. Namely, if the area growth satisfies\footnote{this will be proven in section \ref{eucvsriem}.}
\be\label{qqag}
\mathscr{H}^2_h\big(S\cap B^m_r(p)\big)\;\leq\;\Theta\,r^q\:,\qquad\forall\:\:p\in\vec{\xi}(\di)\:,
\ee
for some universal constant $\Theta$, and some $0<q<2(1+\tau)$, where $\tau$ is as in (\ref{condAF}), then indeed  
\be\label{teague2}
\int_S|\bA^{h_0}_S|_{h_0}^2\,d\mu_{h_0}\;<\;\infty\:
\ee
is true for some universal constant $\Theta$ (which will often be omitted and set to $1$). We will verify that (\ref{teague}) and (\ref{qqag}) together imply (\ref{teague2}). In turn, a classical result by Huber \cite{Hub} (see also \cite{Whi}), guarantees that $S$ is of finite topological type: it is homeomorphic to $\bar{S}\setminus\{a_1,\ldots,a_k\}$, where $\bar{S}$ is a compact surface and $\{a_i\}_{i=1,\dots,k}$ is a set of points. We will be concerned with understanding the surface $S$ around one of these points. For this reason, we suppose there is only such point and we label it $0$. Our surface $S$ may thus be reduced to a connected, oriented, immersed, punctured disk in $\R^m$. The immersion will be denoted by $\bxi:\di\setminus\{0\}\rightarrow(\R^m,h)$. We will suppose that $\bxi$ is a {\it weak immersion} \cite{Riv1}, that is $\bxi$ is Lipschitz and its Gauss map $\bn_{\bxi}$ lies in the Sobolev space $W^{1,2}(D_1(0))$. Moreover, we suppose that
\be\label{condxi}
\bxi(D_r(0))\:\:\:\text{is non-compact}\:\:\:\forall\:\:r\in(0,1)\:,
\ee
and that
\be\label{finiteenergy}
\int_{D_1(0)}|\bA^h_{\bxi}|^2_h\,d\text{vol}_{\bxi^*h}\;<\;\infty\:,
\ee
where $\bA^h_{\bxi}$ is the second fundamental form of $\bxi$. We also impose the  area growth condition:
\be\label{qag}
\mathscr{H}^2_h\big(\bxi(D_1(0))\cap B^m_r(p)\big)\;\leq\;\Theta\,r^q\:,\qquad\forall\:\:p\in\vec{\xi}(\di)\:,
\ee
for some $q<2(1+\tau)$, and $\tau$ as is in (\ref{condAF}).\\ 
A sharpened version of Huber's result due to Stefan M\"uller and Vladimir Sverak \cite{MS} will also be useful to obtain some first information about the asymptotic behavior of the immersion near the branch point located at the origin of the unit disk. More precisely: 
\begin{Prop}\label{introms}
Let $\vec{\xi}:\Om\rightarrow(\R^m,h)$ be a weak immersion into Euclidean space equipped with the asymptotically flat Riemannian metric $h$ satisfying (\ref{condAF}). Suppose that the image of $\vec{\xi}$ is non-compact, complete, has square-integrable fundamental form (\ref{finiteenergy}), area growth (\ref{qag}), and satisfies (\ref{condxi}). Then the immersion is proper, and there exists a reparametrization of the immersion, still denoted $\bxi$, such that $\bxi$ is conformal. Moreover, for an integer $\theta_0\ge1$, it holds
\bes
|\vec{\xi}|_h(x)\;\simeq\;|x|^{-\theta_0}\qquad\text{and}\qquad|\nabla\vec{\xi}|_h(x)\;\simeq\;|x|^{-1-\theta_0}\:,\qquad|x|\ll1\:.
\ees
Here $\nabla$ is the flat gradient with respect to the variable $x$ parametrizing the unit disk. 
\end{Prop} 

\begin{Rm}
As we will see, Proposition \ref{introms} implies that the surface has quadratic area growth:
\bes
\mathscr{H}_h^2\big(\bxi(D_1(0))\cap B^m_R(p)\big)\;\leq\;\Theta R^2\:,
\ees
for some universal constant $\Theta$, and for all radii $R>0$ and all points $p$. 
\end{Rm}

The main object of study of this paper are {\it Willmore surfaces}, which are the critical points of the Willmore energy
\bes
\int_{\Sigma}|\vec{H}^h_{\bxi}|^2_{h}\,d\text{vol}_{\vec{\xi}^*h}\:.
\ees
Clearly, minimal surfaces are Willmore surfaces, so all of our results will in particular apply to complete non-compact minimal surfaces in asymptotically flat space. Being a critical point of the Willmore energy improves the asymptotic behavior of the immersion $\bxi$. 
Imposing only on the ambient metric a general decay to the flat metric as given in the condition (\ref{condAF}), it is possible to show that the second fundamental form of a Willmore surface with finite energy and area growth of type (\ref{qag}) has certain decay properties, as stated in the following theorem, which is our first main result. 
\begin{Th}\label{introeps}
Let the weak Willmore immersion $\vec{\xi}:\Om\rightarrow(\R^m,h)$, the metric $h$, and the integer $\theta_0\ge1$ be as in Proposition \ref{introms}. Then
\be\label{coco}
|\vec{A}^h_{\vec{\xi}}|(p)\;\lesssim\;|p|^{-1}\mu(|p|)\:,\qquad\forall\:\:p\in\vec{\xi}(\di)\quad\text{with}\quad|p|\gg1\:,
\ee
where $\lim_{|p|\nearrow\infty}\mu(|p|)=0$.
\end{Th}
\smallskip
The decay rate given in Theorem \ref{introeps} is unfortunately not sufficient to guarantee that the tangent cone at infinity is unique. In order to reach such a result, as well as for reasons pertaining to applications relevant in general relativity, one must improve (\ref{coco}). To this end, it is necessary to impose further decay on the metric $h$, and demand that it be ``flatter" than the mere (\ref{condAF}). In particular, if we suppose that the decay of the metric $h$ is appropriately synchronized with the asymptotic behavior of $\vec{\xi}$, it is possible to improve (\ref{coco}). This is the content of the next result. 
\begin{Th}\label{introsynchro}
Let the weak Willmore immersion $\vec{\xi}:\Om\rightarrow(\R^m,h)$, the metric $h$, and the integer $\theta_0\ge1$ be as in Proposition \ref{introms}, with the additional assumption that
\be\label{synchro}
h_{\al\beta}(p)\;=\;\delta_{\al\beta}+\text{O}_2(|p|^{-\tau})\qquad\text{for some}\:\:\tau>1-\frac{1}{\theta_0}\quad\text{and for}\:\:\:|p|\gg1\:.
\ee
Then we have for all $\epsilon'>0$:
\bes
|\vec{A}^h_{\vec{\xi}}|(p)\;\lesssim\;|p|^{-1-\frac{1}{\theta_0}+\epsilon'}\:,\qquad\forall\:\:p\in\vec{\xi}(\di)\quad\text{with}\quad|p|\gg1\:.
\ees
Furthermore, in conformal parametrization, $\vec{\xi}$ has near the origin the asymptotic behavior
\be\label{coggy}
\vec{\xi}(x)\;=\;\Re\big(\vec{a}\, x^{-\theta_0}+\vec{a}_1\,x^{1-\theta_0}+\vec{a}_2\,|x|^{-2\theta_0}x^{1+\theta_0}\big)+\text{O}_2\big(|x|^{\theta_0(\tau-1)-\epsilon'}+|x|^{2-\theta_0-\epsilon'}\big)\:,\qquad\forall\:\:\epsilon'>0\:,
\ee
where $\vec{a}$, $\vec{a}_1$, $\vec{a}_2$ are constant vectors in $\C^m$. Here $x$ is to be understood as $x^1+ix^2\in\di$, and $\vec{a}=\vec{a}_R+i\vec{a}_I\in\R^{2}\otimes\C^m$ is a nonzero constant vector satisfying
\be\label{propa}
|\vec{a}_R|_h\,=\,|\vec{a}_I|_h\:\:,\quad
\langle\vec{a}_R,\vec{a}_I\rangle_{h}\,=\,0\:\:,\quad\text{and}\quad \pi_{\bn_h(0)}\vec{a}\,=\,\vec{0}\:.
\ee
Moreover $\pi_{\bn_h(0)}$ denotes the projection onto the normal space of $\vec{\xi}(\di)$ at the point $x=0$. \\[1ex]
Naturally, depending upon the relative sizes of $\tau$ and $\theta_0$, one or more terms in the expansion (\ref{coggy}) are to be absorbed in the most relevant of the two remainders.
\end{Th}

Examples of branched minimal surfaces show that this result is optimal up to the error $\epsilon'>0$. 
A remarkable special case of Theorem \ref{introsynchro} occurs when the surface under study is an {\it embedding}. In that case, it is apparent from the asymptotics given in Proposition \ref{introms} that necessarily $\theta_0=1$, and thus the synchronisation hypothesis (\ref{synchro}) holds for any $\tau>0$.  We feel it is worth rewriting the previous theorem in this special setting. 
\begin{Co}
Let $\vec{\xi}:\Om\rightarrow(\R^m,h)$ be a weak Willmore \underline{embedding} into Euclidean space equipped with the asymptotically flat Riemannian metric $h$ satisfying (\ref{condAF}). Suppose that the image of $\vec{\xi}$ is non-compact, complete, that it has square-integrable fundamental form (\ref{finiteenergy}), area growth (\ref{qag}), and that it satisfies (\ref{condxi}). Then for all $\epsilon'>0$, we have
\bes
|\vec{A}^h_{\vec{\xi}}|(p)\;\lesssim\;|p|^{-2+\epsilon'}\:,\qquad\forall\:\:p\in\vec{\xi}(\di)\quad\text{with}\quad|p|\gg1\:.
\ees
Furthermore, in conformal parametrization, $\vec{\xi}$ has near the origin of the unit disk the asymptotic behavior
\bes
\vec{\xi}(x)\;=\;\Re\big(\vec{a}\, x^{-1}\big)+\text{O}_2(|x|^{-1+\tau})\:,
\ees
where $\vec{a}$ is as in (\ref{propa}).
\end{Co}

Aside from the case $\theta_0=1$, the synchronized hypothesis (\ref{synchro}) might seem somewhat artificial -- although, the authors contend, it is decisive -- for it ties together the asymptotic behavior of the ambient metric $h$ to that of the surface. To obliterate this drawback, it is necessary to assume that the decay of metric $h$ to the Euclidean metric is yet faster, namely we suppose that $h$ is asymptotically Schwarzschild:
\be\label{chinook}
h_{\al\beta}(p)\;=\;\big(1+c|p|^{-1}\big)\delta_{\al\beta}+\text{O}_2(|p|^{-1-\kappa})\qquad\text{for}\:\:\:|p|\gg1\:,
\ee
for some constant $c$ and some $\kappa\in(0,1)$. 
As far the authors know, when $\theta_0\ge2$, it is not possible to significantly improve the asymptotic expansion (\ref{coggy}), even under the stronger hypothesis (\ref{chinook}). However, when $\theta_0=1$, i.e. when the surface is embedded, slightly more can be said. 

\begin{Th}\label{introschwa}
Let the weak Willmore embedding $\vec{\xi}:\Om\rightarrow(\R^m,h)$ be as in Proposition \ref{introms} with $\theta_0=1$, and let the metric $h$ satisfy (\ref{chinook}). \\
Then for all $\epsilon'>0$, near the origin of the unit disk, the conformal parametrization $\vec{\xi}$ has the asymptotic behavior
\bes
\vec{\xi}(x)\;=\;\Re\Big(\vec{a}\, x^{-1}+\vec{a}_1+\vec{a}_2|x|^2x^{-2} \Big)+\vec{c}_0\log|x|^{2}+\text{O}_2(|x|^{\kappa-\epsilon'})\:,
\ees
where $\vec{a}$ is as in Theorem \ref{introsynchro}, while $\vec{a}_1$, $\vec{a}_2$, and $\vec{c}_0$ are constant vectors in $\C^m$. The constant vector $\vec{c}_0$ is normal near the origin:
\be\label{propc}
\pi_{\bn_h(0)}\vec{c}_0\,=\,\vec{c_0}\:.
\ee
\end{Th}

This holds for Willmore immersions and thus in particular for minimal immersions. In the latter case, more can actually be obtained since $\bxi$ is harmonic. Owing to the properties of the vectors $\vec{a}$ and $\vec{c}_0$ given in (\ref{propa}) and (\ref{propc}), one can show that the image of $\vec{\xi}$ can be written as a simple graph over $\R^2\setminus D_R(0)$, for some large enough $R>0$.
 
\begin{Co}\label{intromin}
Let $\vec{\xi}$ be a minimal embedding into Euclidean space $\R^m$ equipped with a Riemannian metric $h$ satisfying the asymptotically Schwarzschild condition (\ref{chinook}). 
Suppose that the image of $\vec{\xi}$ is non-compact, complete, has square-integrable fundamental form (\ref{finiteenergy}), area growth (\ref{qag}), and satisfies (\ref{condxi}). Then for $R$ large enough, the image of $\vec{\xi}$ can be written as a graph over $\R^2\setminus D_R(0)$, namely for all $\epsilon'>0$, it holds:
\bes
(r,\varphi)\,\longmapsto\,\big(r\cos\varphi\,,r\sin\varphi\,,\vec{c}_0\log r+\vec{a}_0+\text{O}_2(r^{\kappa-\epsilon'})\big)\:,
\ees
in the range $\varphi\in[0,2\pi)$ and $r>R$, for some $R$ chosen large enough, and for some $\R^m$-valued constant vectors $\vec{c}_0$ and $\vec{a}_0$.
\end{Co}
With this last statement, we recover Alessandro Carlotto's extension \cite{Car1, Car2}  to complete minimal surfaces in asymptotically Schwarzschild space of Richard Schoen's classical result \cite{Sch} about the end of a complete minimal surface in Euclidean space $\R^3$. Our version is slightly more general as it encompasses minimal surfaces in any codimension. Moreover, whereas in \cite{Car1, Car2}, a \underline{quadratic} area growth condition is also imposed, in our work, this hypothesis is weakened to a $q$-type area growth (\ref{qag}) for any $q<4$.

\subsection{Reformulation of the Problem}\label{reform}

The angle of attack chosen in this paper is as follows. As the metric $h$ is asymptotically flat and our surface satisfies the area growth condition (\ref{qag}), we will first obtain that the immersion $\vec{\xi}$ has square-integrable second fundamental form with respect to the standard Euclidean metric on $\R^m$. A classical result of M\"uller and Sverak \cite{MS} (see also \cite{Hub}) guarantees that $\vec{\xi}$ may be reparametrized into an immersion which is {\it conformal} with respect to the flat metric. For notational convenience, we continue to denote the so-obtained reparametrized immersion by $\vec{\xi}$.  The strategy then consists in ``folding back" the end of the Willmore surface and study the resulting surface, which is the image of an immersion of the punctured unit disk with a singularity at the origin. The main problem in this strategy is to guarantee that the inverted surface satisfies an appropriate variational problem. If the ambient metric were Euclidean, there would be no major problem. Indeed, Willmore surfaces are known to remain Willmore surfaces (possibly singular at a finite set of isolated points) once inverted. This is because inversion in $\R^m$ is a conformal transformation. The presence of the metric $h$ destroys this argument. However, because $h$ is nearly Euclidean in the ``far space", it is possible to apply an inversion to the intersection of $\R^m$ with the complement of a large enough ball. The resulting surface satisfies a perturbed Willmore equation. Using Noether's theorem and its corresponding conservation laws, the Willmore equation, which is {\it a-priori} a fourth-order system, can be recast into a second-order larger system with good analytical dispositions. This technique was originally devised in \cite{Riv1} and made more precise in \cite{Ber2}.

\subsubsection{Euclidean versus Riemannian descriptions}\label{eucvsriem}

We can of course view our immersion $\vec{\xi}$ into $(\R^m,h)$ as an immersion into $(\R^m,h_0)$, where $h_0$ stands for the standard Euclidean metric in $\R^m$. We will respectively denote by $\tilde{h}$ and by $\tilde{h}_0$ the induced metrics $\bxi^*h$ and $\bxi^*h_0$. 
Let us observe once and for all, that 
\bes
|\vec{w}|\;\simeq\;|\vec{w}|_h\qquad\forall\:\:\vec{w}\in \R^m\:.
\ees
We write $a\simeq b$ is to mean that the ratios $|a/b|$ and $|b/a|$ remain bounded as $x$ approaches the origin of $\di$, i.e. as $\bxi(x)$ approaches $\infty$. \\

The goal of this paragraph is to show that the integrability of the second fundamental form $|\vec{A}^h_{\bxi}|^2_h$ along with the hypothesis (\ref{qag}) imply the integrability of $|\vec{A}^{h_0}_{\bxi}|^2$, where $\vec{A}^{h_0}_{\bxi}$ is the second fundamental form of the immersion $\bxi$ into $(\R^m,h_0)$. We begin by inspecting the Gauss maps\footnote{$\star_h$ and $\star$ are the Hodge-star operators associated respectively with the metrics $h$ and $h_0$ in $\R^m$.}:
\bes
\bn_h\;:=\;\star_h\,\dfrac{\partial_{x^1}\bxi\wedge\partial_{x^2}\bxi}{\big|\partial_{x^1}\bxi\wedge\partial_{x^2}\bxi\big|_h}\qquad\text{and}\qquad \bn_{0}\;:=\;\star\,\dfrac{\partial_{x^1}\bxi\wedge\partial_{x^2}\bxi}{\big|\partial_{x^1}\bxi\wedge\partial_{x^2}\bxi\big|}\:.
\ees
One verifies that
\begin{eqnarray}\label{then}
\bn_h&=&\bn_{0}+\bigg(\dfrac{|\tilde{h}|}{|\tilde{h}_0|}-1\bigg)\bn_{0}+|\tilde{h}|^{-1}(\star_h-\star)\big(\partial_{x^1}\bxi\wedge\partial_{x^2}\bxi  \big)\nonumber\\[1ex]
&=&\bn_{0}+\text{O}_2\big(|\bxi|^{-\tau}\big)\:,
\end{eqnarray}
where $\tau$ is as in (\ref{condAF}). \\
For every choice of a $p$-vector $\al$ and a $q$-vector $\beta$ (with $p\ge q$), the interior multiplication $\res_h$ between $\al$ and $\beta$ is implicitly defined through the identity
\bes
\langle \al\res_h \beta,\gamma\rangle_h\;=\;\langle\al,\beta\wedge\gamma\rangle_h\qquad\forall\:\:(p-q)\text{-vector $\gamma$}\:.
\ees
As shown in \cite{MR}, the normal projection of an arbitrary $1$-vector $\bw$ satisfies
\bes
\pi_{\bn_h}\bw\;=\;(-1)^{m-1}\,\bn_h\res_h(\bn_h\res_h\bw)\:.
\ees
The projection $\pi_{\bn_{0}}\bw$ is defined mutatis mutandis, only with respect to the standard Euclidean metric $h_0$ on $\R^m$. With these definitions, it can be verified without much difficulty that for all $\bw$, it holds\footnote{further elaborations in codimension 1 are found in \cite{MSc}.}:
\be\label{projn}
\big|\pi_{\bn_h}\bw-\pi_{\bn_0}\bw\big|\;\lesssim\;\big(|\bn_h|+1\big)\Big(|\bn_h-\bn_0|+\text{O}\big(|\bxi|^{-\tau}\big)|\bn_h|\Big)\,|\bw|\;=\;\text{O}\big(|\bxi|^{-\tau}\big)\,|\bw|\:.
\ee

\smallskip
We let $\nabli{0}{}$, $\nabli{\tilde{h}_0}{}$, $\nabli{h}{}$, and $\nabli{\tilde{h}}{}$ respectively denote the covariant derivatives of the flat Euclidean metric on $\R^2$, and of the metrics $\tilde{h}_0$, $h$, and $\tilde{h}$. The corresponding Christoffel symbols $\gammi{\tilde{h}_0}{}{}$, $\gammi{h}{}{}$, and $\gammi{\tilde{h}}{}{}$ are defined analogously. 
By definition, we have
\bes
\vec{A}^{h_0}_{\bxi}(\partial_{x^i}\bxi,\partial_{x^j}\bxi)\;=\;\nabli{h_0}{\partial_{x^i}\bxi}{\,\partial_{x^j}\bxi}\,-\,\nabli{\tilde{h}_0}{\partial_{x^i}\bxi}{\,\partial_{x^j}\bxi}\;=\;\partial_{x^ix^j}^2\bxi\,-\,\gammi{\tilde{h}_0}{ij}{k}\partial_{x^k}\bxi\:;
\ees
and thus
\begin{eqnarray*}
\vec{A}^{h}_{\bxi}(\partial_{x^i}\bxi,\partial_{x^j}\bxi)&=&\nabli{h}{\partial_{x^i}\bxi}{\,\partial_{x^j}\bxi}\,-\,\nabli{\tilde{h}}{\partial_{x^i}\bxi}{\,\partial_{x^j}\bxi}\;\;=\;\;\partial_{x^ix^j}^2\bxi\,-\,\gammi{h}{\beta\gamma}{\al}\partial_{x^i}\Xi^\beta\partial_{x^j}\Xi^\gamma\vec{E}_\al\,-\,\gammi{\tilde{h}}{ij}{k}\partial_{x^k}\bxi\nonumber\\[1ex]
&=&\vec{A}^{h_0}_{\bxi}(\partial_{x^i}\bxi,\partial_{x^j}\bxi)\,-\,\gammi{h}{\beta\gamma}{\al}\partial_{x^i}\Xi^\beta\partial_{x^j}\Xi^\gamma\vec{E}_\al\,-\,\gammi{\tilde{h}}{ij}{k}\partial_{x^k}\bxi\,+\gammi{\tilde{h}_0}{ij}{k}\partial_{x^k}\bxi\:,
\end{eqnarray*}
where $\Xi^\al$ are the components of $\bxi$ in a fixed orthonormal basis $\{\vec{E}_\al\}_{\al=1,\ldots,m}$ of $\R^m$. Repeated Greek indices indicate summation over $1$ to $m$, while repeated Latin indices indicate summation over $1$ and $2$. For notational convenience, we set $\big(\vec{A}^{h}_{\bxi}\big)_{ij}:=\vec{A}^{h}_{\bxi}(\partial_{x^i}\bxi,\partial_{x^j}\bxi)$.
Projecting the latter on the Euclidean normal space spanned by $\bn_0$ shows that
\be\label{estimo}
\pi_{\bn_0}\big(\vec{A}^{h}_{\bxi}\big)_{ij}\;=\;\big(\vec{A}^{h_0}_{\bxi}\big)_{ij}\,-\,\gammi{h}{\beta\gamma}{\al}\partial_{x^i}\Xi^\beta\partial_{x^j}\Xi^\gamma\pi_{\bn_0}\vec{E}_\al\:.
\ee
The asymptotic form of the metric $h$ given by (\ref{condAF}) implies that
\bes
|\gammi{h}{}{}|\;=\;\text{O}\big(|\bxi|^{-1-\tau}\big)\:.
\ees
It then easily follows from (\ref{estimo}) that
\bes
\big|\vec{A}^{h_0}_{\bxi}\big|^2\;\lesssim\;\big|\vec{A}^{h}_{\bxi}\big|^2_{h}\,+\,|\bxi|^{-2-2\tau}\:,
\ees
where $\nabla\bxi:=\big(\partial_{x^1}\bxi,\partial_{x^2}\bxi\big)$. Using that
\bes
\bigg|\dfrac{|\tilde{h}|^{1/2}}{|\tilde{h}_0|^{1/2}}-1\bigg|\;=\;\text{O}\big(|\bxi|^{-\tau}  \big)\;\ll\;1\:,
\ees
we obtain
\be\label{inter1}
\int_{\bxi(\di)}\big|\vec{A}^{h_0}_{\bxi}\big|^2\,d\text{vol}_{\tilde{h}_0}\;\lesssim\;\int_{\bxi(\di)}\big|\vec{A}^{h}_{\bxi}\big|^2_{h}\,d\text{vol}_{\tilde{h}}\,+\int_{\bxi(\di)}|\bxi|^{-2-2\tau}\,d\text{vol}_{\tilde{h}}\:.
\ee

The first summand on the right-hand side of (\ref{inter1}) is bounded by hypothesis (\ref{finiteenergy}). In light of hypothesis (\ref{qag}), we will now investigate the second summand on the right-hand side of (\ref{inter1}) and verify that it is bounded. 

Using that $\vec{\xi}(\di)\subset\R^m\setminus B_1^m(0)$, we get
\begin{eqnarray*}
\int_{\bxi(\di)}|\bxi|^{-2-2\tau}d\text{vol}_{\bxi^*h}&=&\sum_{j\ge0}\,\int_{\vec{\xi}(\di)\cap (B^m_{2^{j+1}}(0)\setminus B^m_{2^{j}}(0))}|\vec{\xi}|^{-2-2\tau}d\mathscr{H}^2_{h}\nonumber\\[1ex]
&\leq&\sum_{j\ge0}\,2^{-2(1+\tau)j}\mathscr{H}_h^2\big(\vec{\xi}(\di) \cap B^m_{2^{j+1}}(0)\big)\:.
\end{eqnarray*}
The $q$-type area growth given in (\ref{qqag}) then gives
\be\label{charlie007}
\int_{\bxi(\di)}|\bxi|^{-2-2\tau}d\text{vol}_{\bxi^*h}\;\lesssim\;2^q\sum_{j\ge0}2^{(q-2-2\tau)j}\;<\;\infty\:.
\ee
This guarantees that the second summand on the right-hand side of (\ref{inter1}) is bounded, and thus that $\bxi$ has square-integrable second fundamental form as an immersion into the flat Euclidean space $(\R^m,h_0)$. Moreover, by hypothesis, we know that $\bxi$ is complete with $\bxi(D_r(0))$ being non-compact for all $r>0$. We may now call upon the result in \cite{MS} to infer that $\bxi$ may be reparametrized into a {\it proper conformal} immersion of the unit disk into $(\R^m,h_0)$. This reparametrization will simply be denoted $\bxi$, for convenience. Moreover, as shown in \cite{MS}, there exists an integer $\theta_0\ge1$ such that:
\be\label{muls}
|\bxi|(x)\;\simeq\;|x|^{-\theta_0}\qquad\text{and}\qquad|\nabla\bxi|(x)\;\simeq\;|x|^{-\theta_0-1}\:,\:\:\quad|x|\ll1\:,
\ee
where, as before and throughout this paper, $\nabla$ denotes the flat gradient with respect to the variable $x$ parametrizing the unit disk. 

\begin{Rm} From the work of M\"uller and Sverak \cite{MS}, more is known about the conformal factor $\text{e}^{\sigma}$. Namely, 
\bes
\text{e}^{\,\sigma(x)}\;=\;\text{e}^{\,\sigma_0}|x|^{-\theta_0-1}+\text{o}\big(|x|^{-\theta_0-1}\big)\:,
\ees
where $\sigma_0$ is a finite number. Hence, in particular,
\bes
|\bxi|^2(x)\;=\;\text{e}^{\,\sigma_0}|x|^{-\theta_0}+\text{o}\big(|x|^{-\theta_0}\big)\:.
\ees
Let $R>1$ be sufficiently large, and let $r_R$ be such that $r_R^{\theta_0}:=\sigma_0/R$. Note that (\ref{muls}) yields
\begin{eqnarray*}
\int_{\{x\in D_1(0)\,|\,|\bxi|^2(x)\le R\}}|\nabla\bxi|^2(x)\,dx&=&\big(1+\text{o}(1)\big)\,\int_{D_1(0)\setminus D_{r_R}(0)}|\nabla\bxi|^2(x)\,dx\nonumber\\[1ex]
&\simeq&\int_{D_1(0)\setminus D_{r_R}(0)}|x|^{-2\theta_0-2}\,d|x|\;\;\simeq\;\;r_R^{-2\theta_0}\;\;\simeq\;\;R^{2}\:.
\end{eqnarray*}
Since the quantity on the left-hand side of the latter is the area of the surface $\bxi(D_1(0))$ restricted to the ball $B^m_R(0)$, we obtain the quadratic area growth:
\bes
\mathscr{H}^2\big(\bxi(D_1(0))\cap B^m_R(0)\big)\;\lesssim\;R^2\:,
\ees
up to an irrelevant multiplicative constant. Naturally, choosing $R$ large enough and calling upon the fact that the ambient metric $h$ is nearly flat at infinity, we deduce 
\be\label{quadra}
\mathscr{H}_h^2\big(\bxi(D_1(0))\cap B^m_R(0)\big)\;\lesssim\;R^2\:.
\ee
\end{Rm}

\subsubsection{Folding back the surface}

Let $I$ denote the inversion in $\R^m$ about the origin, namely $I(p)=\dfrac{p}{|p|^2}=:y$. One easily verifies that
\bes
g_{\al\beta}(y)\;:=\;|y|^{4}(I_*h)_{\al\beta}(y)\;=\;\delta_{\al\beta}+\text{O}_2(|y|^\tau)\:,\qquad|y|\ll1\:,
\ees
where $\tau$ is as in (\ref{condAF}). \\
We let $\bp:=I\circ\vec{\xi}:\di\rightarrow(B^m_1(0),g)$, $\Sigma:=\bxi(\di)$, and $\Sigma':=\bp(\di)$. It is readily seen that $\bp$ is conformal with respect to the flat metric on $\R^m$, owing to $\vec{\xi}$ being so as well (cf. previous subsection). Also,
\bes
\bp(\di)\:\:\text{has finite area}\qquad\text{and}\qquad\bp(0)\;=\;\vec{0}\:.
\ees
It holds clearly
\be\label{invpar}
\int_{{\Sigma'}}\big|\vec{H}^g_{\bp}\big|^2_{g}\,d\text{vol}_{\vec{\Phi}^*g}\;=\;\int_{\Sigma}\big|\vec{H}^{I^*g}_{\bxi}\big|^2_{I^*g}\,d\text{vol}_{\bxi^*(I^*g)}\;=\;\int_{\Sigma}\big|\vec{H}^{|p|^{-4}h}_{\bxi}\big|^2_{|p|^{-4}h}\,d\text{vol}_{\bxi^*(|p|^{-4}h)}\:.
\ee
It is shown in \cite{Wei} that
\bes
\Lambda(\vec{\zeta},k)\;:=\;\int_{\vec{\zeta}(S)}\Big[\big|\vec{H}^k_{\vec{\zeta}}\big|^2_{k}+\overline{K}^k(T\vec{\zeta})\Big]\,d\text{vol}_{\vec{\zeta}^*k}+\int_{\partial\vec{\zeta}(S)} \kappa_{\vec{\zeta}^*k}\,dS_{{\vec{\zeta}^*k}}
\ees
is an invariant quantity under conformal changes of the metric of $N$. In this generically written expression, $\vec{\zeta}:S\rightarrow(N,k)$ is an immersion of a two-dimensional surface $S$ into a Riemannian manifold equipped with the metric $k$. The sectional curvature of the ambient manifold $(N,k)$ computed on the tangent space of $\vec{\zeta}(S)$ is denoted by $\overline{K}^k(T\vec{\zeta})$, while $\kappa_{\vec{\zeta}^*k}$ is the geodesic curvature of $\partial\vec{\zeta}(S)$, and $dS_{{\vec{\zeta}^*k}}$ is the induced measure on $\partial\vec{\zeta}(S)$. \\
Using (\ref{invpar}), we thus find
\begin{eqnarray*}
\int_{\Sigma}\big|\vec{H}^{h}_{\bxi}\big|^2_{h}\,d\text{vol}_{\bxi^*h} &=&\Lambda(\bxi,h)-\int_{\Sigma}\overline{K}^{h}(T{\bxi})\,d\text{vol}_{\bxi^*h}\,-\,\int_{\partial\Sigma} \kappa_{\bxi^*h}\,dS_{\bxi^*h}\nonumber\\[1ex]
&=&\Lambda(\bxi,|p|^{-4}h)\,-\int_{\Sigma}\overline{K}^{h}(T{\bxi})\,d\text{vol}_{\bxi^*h}\,-\,\int_{\partial\Sigma} \kappa_{\bxi^*h}\,dS_{\bxi^*h}\nonumber\\[1ex]
&=&\int_{\Sigma}\big|\vec{H}^{|p|^{-4}h}_{\bxi}\big|^2_{|p|^{-4}h}\,d\text{vol}_{\bxi^*(|p|^{-4}h)}+\int_{\Sigma}\overline{K}^{|p|^{-4}h}(T{\bxi})\,d\text{vol}_{\bxi^*(|p|^{-4}h)}-\int_{\Sigma}\overline{K}^{h}(T{\bxi})\,d\text{vol}_{\bxi^*h}\nonumber\\
&&\,+\,\int_{\partial\Sigma} \kappa_{\bxi^*(|p|^{-4}h)}\,dS_{\bxi^*(|p|^{-4}h)}\,-\,\int_{\partial\Sigma} \kappa_{\bxi^*h}\,dS_{\bxi^*h}\nonumber\\[1ex]
&=&\int_{{\Sigma'}}\big|\vec{H}^g_{\bp}\big|^2_{g}\,d\text{vol}_{\vec{\Phi}^*g}+\int_{\Sigma}\overline{K}^{|p|^{-4}h}(T{\bxi})\,d\text{vol}_{\bxi^*(|p|^{-4}h)}-\int_{\Sigma}\overline{K}^{h}(T{\bxi})\,d\text{vol}_{\bxi^*h}\nonumber\\
&&\,+\,\int_{\partial\Sigma} \kappa_{\bxi^*(|p|^{-4}h)}\,dS_{\bxi^*(|p|^{-4}h)}\,-\,\int_{\partial\Sigma} \kappa_{\bxi^*h}\,dS_{\bxi^*h}
\end{eqnarray*}
A well-known identity shows that
\bes
\overline{K}^{|p|^{-4}h}(T{\bxi})\,d\text{vol}_{\bxi^*(|p|^{-4}h)}-\overline{K}^{h}(T{\bxi})\,d\text{vol}_{\bxi^*h}\;=\;\big(\Delta_{\bxi^*h}\log|p|^{2}\big)\,d\text{vol}_{\bxi^*h}\:.
\ees
Accordingly, we find
\begin{eqnarray}\label{compH}
\int_{\Sigma'}\big|\vec{H}^{g}_{\bp}\big|^2_{g}\,d\text{vol}_{\bp^*g} &=&\int_{{\Sigma}}\big|\vec{H}^h_{\vec{\xi}}\big|^2_{h}\,d\text{vol}_{\vec{\xi}^*h}-\int_{\Sigma} \Delta_{\bxi^*h}\log|\bxi|^{2}\,d\text{vol}_{\bxi^*h}\nonumber\\
&&\,-\,\int_{\partial\Sigma} \kappa_{\bxi^*(|p|^{-4}h)}\,dS_{\bxi^*(|p|^{-4}h)}\,+\,\int_{\partial\Sigma} \kappa_{\bxi^*h}\,dS_{\bxi^*h}\:.
\end{eqnarray}
The Gauss equation states that
\be\label{gausseq}
\dfrac{1}{2}\int_{\Sigma'}|\vec{A}^g_{\bp}|^2_{g}\,d\text{vol}_{\bp^*g}\;=\;\int_{\Sigma'}|\vec{H}^g_{\bp}|^2_{g}\,d\text{vol}_{\bp^*g}+\Lambda(\bp,g)\:.
\ee
Note that $\Lambda(\bp,g)=\Lambda(\bp,|y|^{-4}g)=\Lambda(I_*\vec{\xi},I_*h)=\Lambda(\vec{\xi},h)$. Combining this to (\ref{compH}) and (\ref{gausseq}) yields that
\begin{eqnarray}\label{compA}
\dfrac{1}{2}\int_{\Sigma'}\big|\vec{A}^{g}_{\bp}\big|^2_{g}\,d\text{vol}_{\bp^*g} &=&\dfrac{1}{2}\int_{{\Sigma}}\big|\vec{A}^h_{\vec{\xi}}\big|^2_{h}\,d\text{vol}_{\vec{\xi}^*h}-\int_{\Sigma} \Delta_{\bxi^*h}\log|\bxi|^{2}\,d\text{vol}_{\bxi^*h}\nonumber\\
&&\,-\,\int_{\partial\Sigma} \kappa_{\bxi^*(|p|^{-4}h)}\,dS_{\bxi^*(|p|^{-4}h)}\,+\,\int_{\partial\Sigma} \kappa_{\bxi^*h}\,dS_{\bxi^*h}\:.
\end{eqnarray}
The last two summands on the right-hand side can be safely assumed to be bounded. This only requires to choose a ``good cut" that isolates the end of the surface. The first summand on the right-hand side is bounded by hypothesis (\ref{finiteenergy}). The second summand is also bounded. To see this, recall that $\Sigma\equiv\bxi(\di)\subset\R^m\setminus B_1^m(0)$, and let $R>0$ be sufficiently large to guarantee that $|\bxi|$ is bounded away from zero. We have
\bes
\bigg|\int_{\Sigma \cap (B^m_R(0)\setminus B^m_1(0))}\Delta_{\bxi^*h}\log|\bxi|^{2}\,d\text{vol}_{\bxi^*h}\bigg|\;=\;\dfrac{2}{R}\, \int_{\Sigma\cap\partial B^m_R(0)}\big|\partial_{\nu_h}|\bxi|\big|\,dS_{\bxi^*h}+2\, \int_{\Sigma\cap\partial B^m_1(0)}\big|\partial_{\nu_h}|\bxi|\big|\,dS_{\bxi^*h}
\ees
where $\nu_h$ is the exterior unit normal with respect to $h$ along $\partial B^m_R(0)$ in the first summand and along $\partial B^m_1(0)$ in the second summand.\\ 
As the induced metric is bounded, we have
\bes
\big|\partial_{\nu_h}|\bxi|\big|_h\;\leq\;\big|\langle\nu_h\,,d|\bxi|\rangle_h\big|\;\leq\;\big|d|\bxi|\big|_h\;\leq\;C_0\:,
\ees
for some constant $C_0$. As the metric $h$ is almost flat (per (\ref{condAF})), it follows that the latter is true with the Euclidean metric in place of $h$. Hence now the estimate
\be\label{maya}
\int_{\Sigma \cap (B^m_R(0)\setminus B^m_1(0))}\Delta_{\bxi^*h}\log|\bxi|^{2}\,d\text{vol}_{\bxi^*h}\;\leq\;\dfrac{2C_0}{R}\mathscr{H}^1_h(\Sigma\cap\partial B^m_R(0))+C_1\:,
\ee
for some constant $C_1$ independent of $R$. 
With the help of the coarea formula and of Fubini's theorem, there exists $\rho=\rho(R)\in(R,2R)$ such that
\bes
\mathscr{H}^1_h(\Sigma\cap\partial B^m_\rho(0))\;\leq\;\dfrac{2}{R}\,\mathscr{H}^2_h(\Sigma\cap B^m_R(0))\;\leq\;C_2R\:,
\ees
for some universal constant $C_2$. Note that the last estimate follows from (\ref{quadra}).  This gives
\bes
\dfrac{1}{\rho}\mathscr{H}^1_h(\Sigma\cap\partial B^m_\rho(0))\;\leq\;C_2\,\dfrac{R}{\rho}\;\leq\;C_2\:.
\ees
This now shows that the right-hand side of (\ref{maya}) is bounded independently of $R$. Letting $R\nearrow\infty$ gives
We now see that
\bes
\int_{\Sigma}\Delta_{\bxi^*h}\log|\bxi|^{2}\,d\text{vol}_{\bxi^*h}\;\equiv\;\lim_{R\rightarrow\infty}\,\int_{\Sigma \cap (B^m_R(0)\setminus B^m_1(0))}\Delta_{\bxi^*h}\log|\bxi|^{2}\,d\text{vol}_{\bxi^*h}\;<\;\infty\:.
\ees

Summarizing our findings, we see that the proper immersion $\bp:\Om\rightarrow(\R^m,g)$ satisfies 
\bes
\bp(\di)\:\:\text{has finite area}\:\:,\qquad \bp(0)=\vec{0}\:\:,\qquad \int_{\bp(\di)}|\vec{A}^g_{\bp}|^2_g\,d\text{vol}_{\bp^*g}<\infty\:.
\ees
Furthermore, we know it is conformal with respect to the Euclidean metric on $\R^m$. \\

Much in the same way that (\ref{compH}) was derived, one derives that 
\begin{eqnarray*}
&&\int_{\Sigma}\big|\vec{H}^{h}_{\bxi}\big|^2_{h}\,d\text{vol}_{\bxi^*h} \nonumber\\[1ex]
&=&
\int_{\Sigma'}\big|\vec{H}^{g}_{\bp}\big|^2_{g}\,d\text{vol}_{\bp^*g} -\int_{\Sigma'} \Delta_{\bp^*g}\log|\bp|^{2}\,d\text{vol}_{\bp^*g}-\int_{\partial\Sigma'} \kappa_{\bp^*(|y|^{-4}g)}\,dS_{\bp^*(|y|^{-4}g)}+\int_{\partial\Sigma'} \kappa_{\bp^*g}\,dS_{\bp^*g}\:.
\end{eqnarray*}
The last three summands are boundary integrals. Because we will only be concerned with local results, we may safely ignore them from our variational analysis. As $\bxi$ is by hypothesis a weak Willmore immersion, it follows that $\bp$ is likewise weak Willmore. Using the results from \cite{MR}, it is not hard to verify that $\bp$ is in fact smooth outside of the origin (where it fails to be Willmore). We will henceforth suppose that
\bes
\bp\,\in\,C^\infty(D_1(0)\setminus\{0\})\cap C^0(D_1(0))\:.
\ees

In passing, a remark which will be of use in the sequel. Observe that
\begin{eqnarray}\label{charlie0}
\int_{\bp(\di)}|\bp|^{-2+2\tau}d\text{vol}_{\bp^*g}&=&\int_{\bxi(\di)}|\bxi|^{-2-2\tau}d\text{vol}_{\bxi^*h}\;<\;\infty\:,
\end{eqnarray}
where we have used (\ref{charlie007}). \\

Singular Willmore immersions in Euclidean space were studied at length in \cite{BR}. The occurrence of the nearly-flat metric $g$ in the present paper will naturally give rise to a perturbed Willmore equation, and our work will consist mainly in showing that this perturbation can be thwarted to produce results akin to those in \cite{BR}. Analyzing a specific class of singular perturbed Willmore immersions (so-called {\it conformally constrained Willmore} immersions) was done in \cite{Ber1}. On the other hand, the Willmore equation in the Riemannian setting was derived in \cite{MR}. Combining the tools and ideas developed in the aforementioned papers is our strategy.

\subsubsection{Quasiconformality of the immersion $\bp$: proof of Proposition \ref{introms}}\label{quasi}
 
We have seen that the immersion is conformal with respect to the Euclidean metric on $\R^m$. In particular, from (\ref{muls}), we have 
\be\label{inca}
|\bp|(x)\;\simeq\;|x|^{\theta_0}\qquad\text{and}\qquad(\tilde{g}_0)_{ij}\;:=\;\partial_{x^i}\bp\cdot\partial_{x^j}\bp\;\simeq\;|x|^{2(\theta_0-1)}\delta_{ij}\:,\qquad|x|\ll1\:,
\ee
where as before $\theta_0\ge1$ is an integer. \\
Because the metric $g$ is only nearly Euclidean, namely
\be
g_{\al\beta}(y)\;=\;\delta_{\al\beta}+\text{O}_2(|y|^{\tau})\:,\qquad|y|\ll1\:,
\ee
for some $\tau>0$, we cannot expect the induced metric $\tilde{g}:=\bp^*g$ to be conformal. At best, it is quasiconformal. The goal of this section is to produce an orthonormal basis of vectors for $\tilde{g}$ and obtain information about it. \\

We start by defining the quantities
\bes
\sigma\;:=\;\dfrac{1}{4}\Big[\tilde{g}_{11}+\tilde{g}_{22}+2\big(\tilde{g}_{11}\tilde{g}_{22}-\tilde{g}^2_{12}\big)^{1/2}  \Big]
\ees
and
\bes
\mu\;:=\;\dfrac{1}{\sigma}\Big[\tilde{g}_{11}-\tilde{g}_{22}+2i\tilde{g}_{12}\Big]\:.
\ees
Setting $z:=x^1+ix^2$ and $\bar{z}:=x^1-ix^2$, one easily verifies that
\bes
\tilde{g}\;=\;\sigma\big|dz+\mu d\bar{z}\big|^2\:.
\ees
Upon letting $w\in\C$ satisfy the Beltrami equation
\bes
\pzb w\;=\;\mu\,\pz w\:,
\ees
we arrive at the conformal representation
\bes
\tilde{g}\;=\;\dfrac{\sigma}{\pz w}|dw|^2\:.
\ees
Note that
\be\label{jeanpetit}
\tilde{g}_{ij}\;=\;\text{e}^{2\nu}\Big(\delta_{ij}+\text{O}_2\big(|\bp|^{\tau}\big)\mathbb{I}_{ij}\Big)\:,
\ee
where $\mathbb{I}_{ij}=1$ for all $i$ and $j$. Here $\nu$ denotes the conformal parameter of the pull-back of the Euclidean metric by $\bp$. In particular, $\text{e}^{\nu}\simeq|x|^{\theta_0-1}$. Hence,
\bes
|x|^{-2(\theta_0-1)}\sigma\;=\;1+\text{O}_2\big(|\bp|^{\tau}\big)\qquad\text{and}\qquad\mu\;=\;\text{O}_2\big(|\bp|^{\tau}\big)\:.
\ees
An exact expression for $\mu$ shall not be necessary for our purposes. Let $\tilde{\mu}:=\mu$ on $\di$ and $\tilde{\mu}:=0$ on $\C\setminus\di$. Consider the Beltrami problem on $\C\,$:
\be\label{bel}
\pzb f\;=\;\tilde{\mu}\,\pz f\:.
\ee
As $|\tilde\mu|<1$, it is known \cite{AIM, Boj} that there exists a solution with $f(0)=0$ and
\be\label{boj1}
\big\Vert \pz f-1\big\Vert_{L^p(\di)}+\big\Vert\pzb f\big\Vert_{L^p(\di)}\;<\;\infty\:,\qquad\forall\:\:p<\infty\:.
\ee
Hence $f-z$ lies in $\bigcap_{p<\infty}W^{1,p}(\di)$. More can be said. Indeed, 
we have
\be\label{fade}
|\nabla\tilde\mu|\;\lesssim\;|\bp|^{-1+\tau}|\nabla\bp|\;\simeq\;|x|^{\tau\theta_0-1}\:,
\ee
where we have used (\ref{inca}). From this it follows that for some $\eta_0>0$, the function $\tilde{\mu}$ lies in $W^{1,2+\eta_0}$. 
Differentiating the Beltrami equation (\ref{bel}) throughout with respect to $z$ yields an equation of the type
\be\label{mucken}
a^{ij}\partial_{x^i}\partial_{x^j}f\;=\;\text{O}\big(|\nabla\tilde\mu||\nabla f|\big)\:,
\ee
with coefficients $a^{11}=1-\tilde{\mu}$, $a^{22}=1+\tilde{\mu}$, and $a^{12}=2i\tilde{\mu}$. As $|\tilde{\mu}|\ll1$, this equation is uniformly elliptic. \\
Owing to (\ref{boj1}) and (\ref{fade}), the right-hand side lies in $L^{2+\eta}$ for some $\eta_0>\eta>0$. Accordingly, we obtain
\bes
\nabla^2f\,\in\,L^{2+\eta}(\di)\:.
\ees
Coupled to (\ref{boj1}), the latter shows that $f$ is invertible on $\di$ and
\be\label{inv}
f^{-1}(w)\;=\;w\,+\,W^{2,2+\eta}(\di)\qquad\text{when}\:\:0<\eta<\eta_0\:,
\ee
with $\eta_0$ as above. \\
The map $\vec{\Psi}(w):=(\bp\circ f^{-1})(w)$ is a continuous immersion of the unit disk, which lies in $W^{1,\infty}\cap W^{2,2}$. By construction, we have that $\Psi(w)$ is conformal with respect to the metric $g$:
\bes
g_{\al\beta}\partial_{u^i}\Phi^\al\partial_{u^j}\Phi^\beta\;=\;\text{e}^{2\la}\delta_{ij}\:,
\ees
where $\la$ is the conformal parameter and $u^1+iu^2:=w$. Since $z=f^{-1}(w)|_{\di}$, we get from (\ref{inv}) and the fact that $g$ is nearly flat that $\la$ and $\nu$ are equivalent near $w=0$ (i.e. near $z=0$). Thus (\ref{inca}) gives now
\bes
|\vec{\Psi}|(u)\;\simeq\;|u|^{\theta_0}\qquad\text{and}\qquad|\nabla\vec{\Psi}|(u)\;\simeq\;|u|^{\theta_0-1}\:,\qquad |u|\ll1\:.
\ees
 In \cite{MS}, it is shown that
\be\label{ms3}
|x||\partial_{x}\nu|(x)\,\in\,L^\infty(\di)\:.
\ee
We will use this inclusion to obtain the following one:
\be\label{ms7}
|u||\partial_u\lambda|(u)\,\in\,\bigcap_{p<\infty}L^{p}(\di)\:,
\ee
which is equivalent to
\be\label{ms888}
|x||\partial_u\lambda|(u)\,\in\,\bigcap_{1<p<\infty}L^{p}(\di)\:.
\ee
Returning to (\ref{jeanpetit}) and differentiating, we see that
\bes
\text{e}^{-2\nu}|\partial_x\tilde{g}|\;\lesssim\;|\partial_x\nu|+|\bp|^{\tau-1}|\nabla\bp|\;\simeq\;|\partial_x\nu|+|x|^{\tau\theta_0-1}\:.
\ees
The Christoffel symbols of the metric $\tilde{g}$ thus satisfy 
\bes
|\gammi{\tilde{g}}{}{}|(x)\;\lesssim\;|\partial_x\nu|(x)+|x|^{\tau\theta_0-1}\:,
\ees
and accordingly we find from (\ref{ms3}) that
\be\label{clefff}
|x||\gammi{\tilde{g}}{}{}|(x)\,\in\,\bigcap_{p<\infty}L^p(\di)\:.
\ee
Weighting the elliptic equation (\ref{mucken}) by $|x|$ and using (\ref{fade}) gives
\bes
|x|a^{ij}\partial_{x^i}f\partial_{x^j}f\;=\;\text{O}(|\nabla f|)\;\in\;L^\infty\:.
\ees
Per the main result in \cite{DST}, we find in particular that 
\bes
|x|\partial_{x}^2 f\,\in\bigcap_{1<p<\infty}L^p(\di)\:,
\ees
which shows that
\be\label{cluff}
|x|\partial_{x}^2 u\,\in\bigcap_{1<p<\infty}L^p(\di)\:.
\ee
The law of transformation for Christoffel symbols is well-known and it is such that
\bes
|\gammi{\tilde{g}}{}{}|(u)\;\lesssim\;|\gammi{\tilde{g}}{}{}|(x)+|\partial_{x}^2u|(x)\:.
\ees
From (\ref{clefff}) and (\ref{cluff}), it now follows that
\bes
|x||\gammi{\tilde{g}}{}{}|(u)\,\in\,\bigcap_{1<p<\infty}L^{p}(\di)\:.
\ees
As $\tilde{g}$ is conformal in the coordinate chart $u$, its Christoffel symbols are given by the partial derivatives of the conformal parameter $\la$. Then (\ref{ms888}) is proved and so is (\ref{ms7}). \\

We will henceforth in this paper only deal with the immersion $\vec{\Psi}$ and the coordinate chart $\{u^1,u^2\}$ on the unit-disk. However, for notational ease, we continue to denote the immersion by the letter $\bp$ and the coordinates by $\{x^1,x^2\}$. This should not generate any confusion for the reader. From the way it was constructed, it is clear that $\bp$ is a conformal immersion for the metric $g$, and that it lies in the space $W^{2,2}\cap W^{1,\infty}$. It is continuous at the origin with $\bp(0)=\vec{0}$. Its conformal factor $\text{e}^\la$ is comparable to $|x|^{\theta_0-1}$. Its second fundamental form (understood with respect to the metric $g$ or to the flat metric) is bounded in $L^2$. Of course, because this ``new" immersion is merely a reparametrized version of its ``old" self, it continues to be a critical point of the Willmore energy. \\

We shall not prove directly Theorem \ref{introeps}, Theorem \ref{introsynchro}, and Theorem \ref{introschwa} as they are stated in Section I.1. We will instead prove the following counterpart versions, from which the statements given in Section I.1 easily ensue. We will suppose that the ambient metric $g$ satisfies
\be\label{methyp}
g_{\al\beta}(y)\;=\;\delta_{\al\beta}+\text{O}_2(|y|^\tau)\:,\qquad|y|\ll1\:,
\ee 
for some $\tau>0$ as in (\ref{condAF}). \\
The conformal immersion $\bp$ is a critical point of the Willmore functional and it satisfies
\be\label{resy1}
\bp\in C^\infty(\di\setminus\{0\})\cap C^0(\di)\:,\quad \bp(0)=\vec{0}\:,\quad \int_{\di}|\vec{A}^g_{\bp}|^2_g\,d\text{vol}_{\bp^*g}<\infty\:.
\ee
Moreover, its conformal parameter satisfies for some integer $\theta_0\ge1$:
\be\label{resy2}
\text{e}^{\la}(x)\simeq|x|^{\theta_0-1}\quad\text{and}\quad |x|\nabla\la\in\bigcap_{p<\infty}L^p(\di)\:.
\ee
Repeating mutatis mutandis (\ref{inter1}) and with the help of (\ref{charlie0}), one easily deduces that
\be\label{tango1}
\int_{\di}|\vec{A}^{h_0}_{\bp}|^2\,d\text{vol}_{\bp^*h_0}<\infty\:,
\ee
where, as before, $h_0$ stands for the usual Euclidean metric on $\R^m$. 

\begin{Th}\label{eps}
Let the conformal Willmore immersion ${\bp}:\Om\rightarrow(\R^m,g)$, the metric $g$, and the integer $\theta_0\ge1$ be as in (\ref{methyp})-(\ref{resy2}). Then
\bes
|\vec{A}^g_{\bp}|(y)\;\lesssim\;|y|^{-1}\mu(|y|)\:,\qquad\forall\:\:y\in\bp(\di)\quad\text{with}\quad|y|\ll1\:,
\ees
where $\lim_{|y|\searrow0}\mu(|y|)=0$.
\end{Th}

\begin{Th}\label{synchroth}
Let the Willmore conformal immersion $\bp:\Om\rightarrow(\R^m,g)$, the metric $g$, and the integer $\theta_0\ge1$ be as in (\ref{methyp})-(\ref{resy2}), with the additional assumption that
\bes
g_{\al\beta}(y)\;=\;\delta_{\al\beta}+\text{O}_2\big(|y|^{\tau}\big)\qquad\text{for some}\:\:\tau>1-\frac{1}{\theta_0}\quad\text{and for}\:\:\:|y|\ll1\:.
\ees
Then for all $\epsilon'>0$, we have
\bes
|\vec{A}^g_{\bp}|(y)\;\lesssim\;|y|^{-1+\frac{1}{\theta_0}-\epsilon'}\:,\qquad\forall\:\:y\in\bp(\di)\quad\text{with}\quad|y|\ll1\:.
\ees
Furthermore, in parametrization, $\bp$ has near the origin the asymptotic behavior
\bes
\bp(x)\;=\;\Re\big(\vec{B}\, x^{\theta_0}+\vec{B}_1\,x^{\theta_0+1}+\vec{B}_2|x|^{2\theta_0}x^{1-\theta_0}\big)+\text{O}_2\big(|x|^{\theta_0(\tau+1)-\epsilon'}+|x|^{\theta_0+2-\epsilon'}\big)\:,\qquad\forall\:\:\epsilon'>0\:,
\ees
where $\vec{B}$, $\vec{B}_1$, and $\vec{B}_2$ are constant vectors in $\C^m$. Here, $x$ is to be understood as $x^1+ix^2\in\C$, and $\vec{B}=\bB_R+i\bB_I\in\R^{2m}$ is a nonzero constant vector satisfying
\bes
|\vec{B}_R|_g\,=\,|\vec{B}_I|_g\:\:,\quad
\langle\vec{B}_R,\vec{B}_I\rangle_{g}\,=\,0\:\:,\quad\text{and}\quad \pi_{\bn_g(0)}\vec{B}\,=\,\vec{0}\:.
\ees
Moreover, $\pi_{\bn_g(0)}$ denotes the projection onto the normal space of $\bp(\di)$ at the point $x=0$.\\[1ex]
The mean curvature vector has the expansion
\be\label{coggy2}
\bH^g_{\bp}(x)\;=\;-2\vec{\gamma}_0\log|x|+\Re\big(\vec{E}_0x^{1-\theta_0}\big)+\text{O}\big(|x|^{\text{min}\{\tau,2-\theta_0\}-\epsilon'}\big)\qquad\forall\:\:\epsilon'>0\:,
\ee
where $\vec{\gamma}_0\in\R^m$ and $\vec{E}_0\in\C^m$ are constant vectors.\\[1ex]
Naturally, depending upon the relative sizes of $\theta_0$ and $\tau$, one or more summands in the expansion (\ref{coggy2}) are to be absorbed in the remainder. 
\end{Th}

Finally, the pendant of Theorem \ref{introschwa}, namely embeddings (i.e. $\theta_0=1$) in asymptotically Schwarzschild spaces, reads:

\begin{Th}\label{schwat}
Let the Willmore conformal \underline{embedding} $\bp:\Om\rightarrow(\R^m,g)$ satisfy (\ref{resy1}) and let the metric $g$ be such that
\bes
g_{\al\beta}(y)\;=\;\big(1+c|y|\big)\delta_{ij}+\text{O}_2(|y|^{1+\kappa})\:,\qquad|y|\ll1\:,
\ees
for some $0<\kappa<1$ and some constant $c$. 
Then for all $\epsilon'>0$, we have
\bes
|\vec{A}^g_{\bp}|(y)\;\lesssim\;|y|^{\kappa-\epsilon'}\:,\qquad\forall\:\:y\in\bp(\di)\quad\text{with}\quad|y|\ll1\:.
\ees
Furthermore, in parametrization, $\bp$ has near the origin the asymptotic behavior
\bes
\bp(x)\;=\;\Re\Big(\vec{B}x+\vec{B}_1x^2\Big)+\vec{C}_0|x|^{2}\big(\log|x|^{2}-C_1\big)+\text{O}_2\big(|x|^{\kappa+2-\epsilon'}\big)\:,
\ees
where $\vec{B}$ is as in Theorem \ref{synchroth}, while $\bB_1\in\C^m$, $C_1\in\R$ are constant, and $\vec{C}_0$ is a constant vector in $\C^m$ with
\bes
\pi_{T_g(0)}\vec{C}_0\,=\,\vec{0}\:.
\ees
\end{Th}

\noindent
\textbf{Acknowledgments:\:\:} This work was completed during the first author's stay at the Forschungsinstitut f\"ur Mathematik (FIM) of the ETH in Z\"urich. He expresses his gratitude for the excellent working conditions and facilities of the FIM. The authors also wish to thank Alessandro Carlotto for bringing to their attention some imprecisions and ambiguities in an early version of this work, as well as for stimulating fruitful discussions. 

\setcounter{equation}{0} 
\reset

\section{Proofs of the Theorems}

\subsection{The Willmore Equation}

As is shown in \cite{MR}, a conformal immersion $\bp$ (with conformal factor $\la$) in a Riemannian space $(\R^m,g)$, which is a critical point of the Willmore energy $\int|\bH^g_{\bp}|^2_gd\text{vol}_{\bp^*g}$ satisfies the following partial differential equation:
\be\label{www}
\D{g}{}{*}\big(\D{g}{}{}\vec{H}^{g}_{\bp}-2\pi_{\bn_g} \D{g}{}{}\vec{H}^{g}_{\bp}+|\vec{H}^{g}_{\bp}|^2_g\nabla\bp  \big)-\text{e}^{2\la}\big(\widetilde{R}(\vec{H}^{g}_{\bp})-R^\perp_{\bp}(T\bp)\big)\;=\;\vec{0}\:,
\ee
where $\D{g}{}{}$ and $\D{g}{}{*}$ are respectively the covariant gradient and divergence corresponding to the metric $g$, namely
\bes
\D{g}{}{}\vec{f}\;:=\;\big(\nabli{g}{\partial_{x^1}\bp}\vec{f}\,,\nabli{g}{\partial_{x^2}\bp}\vec{f}\big)\qquad\text{and}\qquad \D{g}{}{*}(\vec{u}\,,\vec{v})\;:=\;\nabli{g}{\partial_{x^1}\bp}\vec{u}+\nabli{g}{\partial_{x^2}\bp}\vec{v}\:.
\ees
As before, $\nabli{g}{}$ is the covariant derivative associated with the metric $g$, while $\nabla$ stands for the flat gradient: $\nabla\vec{f}:=(\partial_{x^1}\vec{f},\partial_{x^2}\vec{f})$. The other two terms appearing in the variation of the Willmore energy are defined as follows. 
\be\label{curv00}
\left\{\begin{array}{lcl}
\text{e}^{2\la}\widetilde{R}(\vec{H}^{g}_{\bp})\;=\;-\,\pi_{\bn_g}\Big[\sum_{j=1,2}\text{Riem}^g(\vec{H}^{g}_{\bp},\partial_{x^j}\bp)\partial_{x^j}\bp\Big] \\[1ex]
\text{e}^{2\la}R^\perp_{\bp}(T\bp)\;=\;\Big[\pi_{T_g}\big(\text{Riem}^g(\partial_{x^1}\bp,\partial_{x^2}\bp)\vec{H}^{g}_{\bp}\big)   \Big]^\perp\:,
\end{array}\right.
\ee
where $\pi_{\bn_g}$ and $\pi_{T_g}$ denote respectively the projection onto the normal and onto the tangent space of $\bp$. The operator $\perp$ is intrinsically defined as
\bes
\vec{X}^\perp\;:=\;(\bp_*)\circ \star_g \circ (\bp_*)^{-1}(\bX)\qquad\text{for}\:\:\:\vec{X}\in\bp_*(T\di)\:,
\ees
where $\bp_*$ is the push-forward of $\bp$, and $\star_g$ is the Hodge-star operator corresponding to the metric $g$. \\

Naturally, to us, the equation (\ref{www}) will only hold on $\Om$. The goal will be to understand how $\bH^g_{\bp}(x)$ and $\bp(x)$ behave near the origin $x=0$. 

\subsubsection{The asymptotically flat case and proof of Theorem \ref{eps}}

As before, we suppose that the metric $g$ satisfies (\ref{methyp}). The components of the Riemann tensor of the metric $g$ computed on the surface parametrized by $\bp$ satisfy
\bes
\text{Riem}^g(\vec{u},\vec{v})\vec{w}\;=\;\text{O}\big(|\bp|^{-2+\tau}|\vec{u}||\vec{v}||\vec{w}|\big)\qquad\forall\:\vec{u},\,\vec{v},\,\vec{w}\:.
\ees
Hence (\ref{curv00}) and (\ref{www}) give
\be\label{www1}
\D{g}{}{*}\big(\D{g}{}{}\vec{H}^{g}_{\bp}-2\pi_{\bn_g} \D{g}{}{}\vec{H}^{g}_{\bp}+|\vec{H}^{g}_{\bp}|^2_g\nabla\bp  \big)\;=\;\text{O}\big(|\bp|^{-2+\tau}|\nabla\bp|^2|\vec{H}^g_{\bp}| \big)\qquad\text{on}\:\:D^2\setminus\{0\}\:.
\ee
Using once more the hypothesis on the metric $g$, we also verify that
\begin{eqnarray*}
\D{g}{\partial_{x^j}\bp}{}\vec{f}&:=&\partial_{x^j}\vec{f}+\gammi{g}{\beta\gamma}{\al}\partial_{x^j}\Phi^\beta f^\gamma\vec{E}_\al\;\;=\;\;\partial_{x^j}\vec{f}\,+\text{O}\big(|\bp|^{-1+\tau}|\nabla\bp||\vec{f}|\big)
\end{eqnarray*}
holds for all $\vec{f}$. In this expression, $\gammi{g}{\beta\gamma}{\al}$ are the Christoffel symbols of the metric $g$, while $\Phi^\beta$ and $f^\gamma$ are respectively the components of $\bp$ and of $\vec{f}$ in a fixed basis $\{\vec{E}_\al\}_{\al=1,\ldots,m}$ of $\R^m$. Introducing this information into (\ref{www1}) gives now the following equation holding on the punctured unit disk:
\be\label{equa}
\text{div}\big(\nabla\vec{H}^{g}_{\bp}-2\pi_{\bn_g} \nabla\vec{H}^{g}_{\bp}+|\vec{H}^{g}_{\bp}|^2\nabla\bp+\bw_1\big)\;=\;\bw_2-\vec{E}_\al\sum_{j=1,2}\gammi{g}{\beta\gamma}{\al}\partial_{x^j}\Phi^\beta\big(\partial_{x^j}\vec{H}^g_{\bp}-2\pi_{T_g}\partial_{x^j}\vec{H}^g_{\bp}\big)^\gamma\:,
\ee
where
\be\label{hyppo}
\left\{\begin{array}{lcl}
\bw_1&=&\text{O}\big(|\bp|^{-1+\tau}|\nabla\bp||\vec{H}^g_{\bp}|\big) \\[1.5ex]
\bw_2&=&\text{O}\big(|\bp|^{-1+\tau}|\nabla\bp|^2|\vec{H}^{g}_{\bp}|^2+|\bp|^{-2+\tau}|\nabla\bp|^2|\vec{H}^{g}_{\bp}| \big)\:.
\end{array}\right.
\ee
Note that we have used the simple fact that $\pi_{T_g}=\text{id}-\pi_{\bn_g}$.\\

One checks that
\bes
\pi_{T_g}\D{g}{\partial_{x^j}\bp}{}\bH^g_{\bp}\;=\;-\,\sum_{k=1,2}\big\langle\bH^g_{\bp}\,,\big(\vec{A}^g_{\bp}\big)_{jk}\big\rangle_g\,\partial_{k}\bp\:,
\ees
whence
\be\label{goya1}
\big|\pi_{T_g}\partial_{x^j}\bH^g_{\bp}\big|\;\lesssim\;|\nabla\bp||\bH^g_{\bp}||\vec{A}^g_{\bp}|+|\gammi{g}{}{}||\nabla\bp||\bH^g_{\bp}|\;=\;\text{O}\big(|\nabla\bp||\bH^g_{\bp}||\vec{A}^g_{\bp}|+|\bp|^{-1+\tau}|\nabla\bp||\bH^g_{\bp}|\big)\:.
\ee
On the other hand, we have
\be\label{goya777}
\sum_{j=1,2}\gammi{g}{\beta\gamma}{\al}\partial_{x^j}\Phi^\beta\big(\partial_{x^j}\vec{H}^g_{\bp}\big)^{\!\gamma}\vec{E}_\al\;=\;\text{div}\big(\gammi{g}{\beta\gamma}{\al}\nabla\Phi^\beta(\vec{H}^g_{\bp})^{\!\gamma}\vec{E}_\al \big)\,-\,\big(\nabla\gammi{g}{\beta\gamma}{\al}\cdot\nabla\Phi^\beta+\gammi{g}{\beta\gamma}{\al}\Delta\Phi^\beta  \big)(\bH^g_{\bp})^\gamma\vec{E}_\al\:.
\ee
As was done in section \ref{eucvsriem}, we have
\bes
\vec{A}^{g}_{\bp}(\partial_{x^i}\bp,\partial_{x^j}\bp)\;\;=\;\;\partial_{x^ix^j}^2\bp\,-\,\gammi{g}{\beta\gamma}{\al}\partial_{x^i}\Phi^\beta\partial_{x^j}\Phi^\gamma\vec{E}_\al\,-\,\gammi{\tilde{g}}{ij}{k}\partial_{x^k}\bp\:.
\ees
Since $\tilde{g}_{ij}=\text{e}^{2\la}\delta_{ij}$, we contract this identity and use the well-known fact that for a conformal metric $\tilde{g}^{ij}\big(\gammi{\tilde{g}}{ij}{k}\big)=0$, to find
\be\label{laseulefaute}
2\text{e}^{2\la}\vec{H}^{g}_{\bp}\;=\;\Delta\bp+\text{O}\big(|\gammi{g}{}{}||\nabla\bp|^2\big)\;=\;\Delta\bp\,+\,\text{O}\big(|\bp|^{\tau-1}|\nabla\bp|^2\big)\:,
\ee
where $\Delta$ is simply the flat Laplace operator, and we have used the previously encountered fact that $\gammi{g}{}{}=\text{O}(|\bp|^{\tau-1})$. Brought into (\ref{goya777}), this information yields
\begin{eqnarray*}
&&\bigg|\sum_{j=1,2}\gammi{g}{\beta\gamma}{\al}\partial_{x^j}\Phi^\beta\big(\partial_{x^j}\vec{H}^g_{\bp}\big)^{\!\gamma}\vec{E}_\al\,-\,\text{div}\big(\gammi{g}{\beta\gamma}{\al}\nabla\Phi^\beta(\vec{H}^g_{\bp})^{\!\gamma}\vec{E}_\al \big)  \bigg|\\[-0ex]
&&\hspace{3cm}\;=\;\text{O}\big(|\bp|^{-2+\tau}|\nabla\bp|^2|\bH^g_{\bp}|+|\bp|^{-1+\tau}|\nabla\bp|^2|\bH^g_{\bp}|^2 \big)\:,
\end{eqnarray*}
where we have used that $|\nabla\gammi{g}{}{}|=\text{O}(|\bp|^{-2+\tau}|\nabla\bp|)$. 
Introducing the latter and (\ref{goya1}) into (\ref{equa})-(\ref{hyppo}) gives the following equation which holds on the punctured unit disk:
 \bes
\text{div}\big(\nabla\vec{H}^{g}_{\bp}-2\pi_{\bn_g} \nabla\vec{H}^{g}_{\bp}+|\vec{H}^{g}_{\bp}|^2\nabla\bp+\bu_1\big)\;=\;\bu_2\:,
\ees
where
\bes
\left\{\begin{array}{lcl}
\bu_1&=&\text{O}\big(|\bp|^{-1+\tau}|\nabla\bp||\vec{H}^g_{\bp}|\big) \\[1.5ex]
\bu_2&=&\text{O}\big(|\bp|^{-1+\tau}|\nabla\bp|^2|\vec{H}^g_{\bp}|^2+|\bp|^{-2+\tau}|\nabla\bp|^2|\vec{H}^{g}_{\bp}|\big)\:.
\end{array}\right.
\ees
Since, as seen in section \ref{quasi}, it holds near the origin that
\bes
|\bp|(x)\;\simeq\;|x|^{\theta_0}\qquad\text{and}\qquad |\nabla\bp|(x)\;\simeq\;|x|^{\theta_0-1}\:,
\ees
the problem which we are considering is
 \be\label{equa2}
\text{div}\big(\nabla\vec{H}^{g}_{\bp}-2\pi_{\bn_g} \nabla\vec{H}^{g}_{\bp}+|\vec{H}^{g}_{\bp}|^2\nabla\bp+\bu_1\big)\;=\;\bu_2\qquad\text{on}\:\:D_1(0)\setminus\{0\}\:,
\ee
with
\be\label{hypo27}
\left\{\begin{array}{lcl}
\bu_1&=&\text{O}\big(|x|^{\theta_0(\tau-1)}|\nabla\bp||\vec{H}^g_{\bp}|\big) \\[1.5ex]
\bu_2&=&\text{O}\big(|x|^{\theta_0(\tau-1)}|\nabla\bp|^2|\vec{H}^g_{\bp}|^2+|x|^{\theta_0(\tau-2)}|\nabla\bp|^2|\vec{H}^{g}_{\bp}|\big)\:.
\end{array}\right.
\ee

{\it Mutatis mutandis} (\ref{then}), we have
\bes
\nabla\bn_g\;=\;\nabla\bn_0+\text{O}\big(|\bp|^{-1+\tau}|\nabla\bp|\big)\:.
\ees
We have already seen in (\ref{charlie0}) that $|\bp|^{-1+\tau}|\nabla\bp|$ belongs to $L^2(\di)$. In addition, using (\ref{tango1}), one easily checks that
\bes
\int_{D_1(0)}|\nabla\bn_0|^2dx\;=\;\int_{D_1(0)}|\vec{A}^{h_0}_{\bp}|^2d\text{vol}_{\bp^*h_0}\;<\;\infty\:,
\ees
it follows that $\nabla\bn_g\in L^2(D_1(0))$. Recall that $h_0$ stands for the standard Euclidean metric on $\R^m$.\\
Thus, owing to (\ref{then}), there holds easily
\bes
\nabla\big(\star_g\!\bn_g-\star\bn_0\big)\;=\;\text{O}\big(|\bp|^{\tau}|\nabla\bn_g|+|\bp|^{-1+\tau}|\nabla\bp||\bn_g|  \big)\,\in\,L^2(\di)\:.
\ees
Hence
\bes
\Vert\nabla(\star_g\bn_g)\Vert_{L^2(D_1(0))}\;<\;\infty\:.
\ees

We are only interested in local results around the origin of the punctured disk. Rescaling the domain if necessary, we may and will assume that for some $\eps_0>0$ chosen as small as we deem useful, it holds
\be\label{epshyp}
\int_{D_1(0)}|\nabla(\star_g\bn_g)|^2dx\,+\int_{D_1(0)}|\nabla\bn_g|^2dx\;<\;\eps_0\:,
\ee
without any loss on the quantitative equation (\ref{equa2}) or on the qualitative hypotheses (\ref{hypo27}). We will now prove Theorem \ref{eps}. 
\begin{Lma}\label{epsss}
There holds
\bes
\lim_{r\rightarrow0}\delta(r)\,=\,0\qquad\text{and}\qquad\int_{0}^{1/2}\delta^2(r)\dfrac{dr}{r}\,<\,\infty\:.
\ees
\end{Lma}
{\bf Proof.} The argument relies on a so-called $\eps$-regularity estimate for equations of the type (\ref{equa2}) under the hypothesis (\ref{epshyp}). The formulation found in \cite{BWW} states that
\be\label{epsreg}
\Vert\text{e}^\la\vec{A}^g_{\bp}\Vert_{L^\infty(D_s)}\;\le\;C_0\bigg[s\Vert\text{e}^\la\vec{u}_2\Vert_{L^2(D_{2s})}+s^{1/2}\Vert\text{e}^\la\vec{u}_1\Vert_{L^4(D_{2s})}+\dfrac{1}{s}\,\Vert\nabla\bn_g\Vert_{L^2(D_{2s})}\bigg]\:,
\ee
holds for any flat disk $D_{2s}\subset D_1(0)\setminus\{0\}$, where $C_0$ is a universal constant. \\

Let $r\in(0,1/2)$. Clearly, there exists a finite number of points $x_j\in\partial D_r(0)$ and a positive constant $c<1/4$ such that
\bes
\partial D_r(0)\,\subset\,\bigcup_{j=1}^{N}D_{cr}(x_j)\qquad\text{and}\qquad D_{2cr}(x_j)\,\subset\,D_{2r}(0)\setminus D_{r/2}(0)\:.
\ees
For some point $x_j\in\partial D_r(0)$, we have
\begin{eqnarray}\label{newx1}
r^{\frac{1}{2}}\Vert\text{e}^{\la}\bu_1\Vert_{L^{4}(D_{2cr}(x_j))}&\lesssim&r^{\theta_0\tau-1+\frac{1}{2}}\big\Vert\text{e}^{\la}\bH^g_{\bp}\big\Vert_{L^{4}(D_{2cr}(x_j))}\nonumber\\[1ex]
&\lesssim&r^{\theta_0\tau-1}\big\Vert|x|\text{e}^{\la}\bH^g_{\bp}\big\Vert^{\frac{1}{2}}_{L^\infty(D_{2cr}(x_j))}\Vert\nabla\bn_g\Vert_{L^2(D_{2cr}(x_j))}^{\frac{1}{2}}\nonumber\\[1ex]
&\lesssim&r^{\theta_0\tau-1}\big\Vert|x|\text{e}^{\la}\bH^g_{\bp}\big\Vert^{\frac{1}{2}}_{L^\infty(D_{2r}(0))}\Vert\nabla\bn_g\Vert_{L^2(D_{2r}(0)\setminus D_{r/2}(0))}^{\frac{1}{2}}\:.
\end{eqnarray}
On the other hand, using (\ref{hypo27}), we find
\begin{eqnarray}\label{newx2}
r\Vert\text{e}^{\la}\bu_2\Vert_{L^{2}(D_{2cr}(x_j))}&\lesssim&r^{\theta_0\tau}\big\Vert\text{e}^{\la}\bH^g_{\bp}\big\Vert^{2}_{L^{4}(D_{2cr}(x_j))}+r^{\theta_0\tau}\big\Vert\text{e}^{\la}\bH^g_{\bp}\big\Vert_{L^{2}(D_{2cr}(x_j))}\nonumber\\[1ex]
&\lesssim&r^{\theta_0\tau}\big\Vert|x|\text{e}^{\la}\bH^g_{\bp}\big\Vert_{L^\infty(D_{2r}(0))}\Vert\nabla\bn_g\Vert_{L^2(D_{2r}(0)\setminus D_{r/2}(0))}\nonumber\\
&&\hspace{2cm}+\:r^{\theta_0\tau}\Vert\nabla\bn_g\Vert_{L^2(D_{2r}(0)\setminus D_{r/2}(0))}\:.
\end{eqnarray}
As all quantities involved are assumed to be smooth away from the singularity, we can invoke the estimate (\ref{epsreg}) and use (\ref{newx1}) and (\ref{newx2}) to find
\begin{eqnarray}\label{ouatev}
&&\hspace{-4cm}\Vert\text{e}^{\la}\bA^g_{\bp}\Vert_{L^{\infty}(\partial D_{r}(0))}\;\;\leq\;\;\Vert\text{e}^{\la}\bA^g_{\bp}\Vert_{L^{\infty}(D_{cr}(x_j))}\nonumber\\[1ex]
&&\hspace{-0.9cm}\lesssim\:\;  r^{\theta_0\tau-1}\big\Vert|x|\text{e}^{\la}\bH^g_{\bp}\big\Vert^{\frac{1}{2}}_{L^\infty(D_{2r}(0))}\Vert\nabla\bn_g\Vert_{L^2(D_{2r}(0)\setminus D_{r/2}(0))}^{\frac{1}{2}}\nonumber\\[1ex]
&&\hspace{0cm}+\:\:r^{\theta_0\tau}\big\Vert|x|\text{e}^{\la}\bH^g_{\bp}\big\Vert_{L^\infty(D_{2r}(0))}\Vert\nabla\bn_g\Vert_{L^2(D_{2r}(0)\setminus D_{r/2}(0))}\nonumber\\[1ex]
&&\hspace{1cm}+\:\:r^{\theta_0\tau}\Vert\nabla\bn_g\Vert_{L^2(D_{2r}(0)\setminus D_{r/2}(0))}\nonumber\\[1ex]
&&\hspace{2cm}+\:\:r^{-1}\Vert\nabla\bn_g\Vert_{L^2(D_{2r}(0)\setminus D_{r/2}(0))}\:.
\end{eqnarray}
When $\rho\in(0,1/2)$, we can always find $r\leq\rho$ such that
\bes
\big\Vert|x|\text{e}^{\la}\bA^g_{\bp}\big\Vert_{L^\infty(D_\rho(0))}\;\leq\;r\,\Vert\text{e}^{\la}\bA^g_{\bp}\Vert_{L^\infty(\partial D_r(0))}\:.
\ees
Combining this to (\ref{ouatev}) yields
\begin{eqnarray*}
\big\Vert|x|\text{e}^{\la}\bA^g_{\bp}\big\Vert_{L^\infty(D_\rho(0))}&\lesssim&\rho^{\theta_0\tau}\big\Vert|x|\text{e}^{\la}\bH^g_{\bp}\big\Vert^{\frac{1}{2}}_{L^\infty(D_{2\rho}(0))}\Vert\nabla\bn_g\Vert_{L^2(D_{2\rho}(0)\setminus D_{\rho/2}(0))}^{\frac{1}{2}}\nonumber\\[1ex]
&&\hspace{0cm}+\:\:\rho^{\theta_0\tau+1}\big\Vert|x|\text{e}^{\la}\bH^g_{\bp}\big\Vert_{L^\infty(D_{2\rho}(0))}\Vert\nabla\bn_g\Vert_{L^2(D_{2\rho}(0)\setminus D_{\rho/2}(0))}\nonumber\\[1ex]
&&\hspace{.75cm}+\:\:\rho^{\theta_0\tau}\Vert\nabla\bn_g\Vert_{L^2(D_{2\rho}(0)\setminus D_{\rho/2}(0))}\nonumber\\[1ex]
&&\hspace{1.5cm}+\:\:\Vert\nabla\bn_g\Vert_{L^2(D_{2\rho}(0)\setminus D_{\rho/2}(0))}\:.
\end{eqnarray*}
For notational convenience, we set
\be\label{deldef}
\delta(\rho)\;:=\;\big\Vert|x|\text{e}^{\la}\bA^g_{\bp}\big\Vert_{L^\infty(D_\rho(0))}\:,
\ee
so as to recast the latter in the form
\begin{eqnarray}\label{rebelotte}
\delta(\rho)&\lesssim&\rho^{\theta_0\tau}\Vert\nabla\bn_g\Vert_{L^2(D_{2\rho}(0)\setminus D_{\rho/2}(0))}^{\frac{1}{2}}\,\delta(2\rho)^{\frac{1}{2}}+\rho^{\theta_0\tau+1}\Vert\nabla\bn_g\Vert_{L^2(D_{2\rho}(0)\setminus D_{\rho/2}(0))}\,\delta(2\rho)\nonumber\\[1ex]
&&\hspace{1cm}+\:\;\Vert\nabla\bn_g\Vert_{L^2(D_{2\rho}(0)\setminus D_{\rho/2}(0))}\nonumber\\[1ex]
&\lesssim&\rho^{\theta_0\tau}\big(1+\rho\,\Vert\nabla\bn_g\Vert_{L^2(D_{2\rho}(0)\setminus D_{\rho/2}(0))}\big)\,\delta(2\rho)+\Vert\nabla\bn_g\Vert_{L^2(D_{2\rho}(0)\setminus D_{\rho/2}(0))}\:.
\end{eqnarray}
Using the uniform bound (\ref{epshyp}) and the fact that $\theta_0\tau>0$ gives first that
\bes
\delta(\rho)\;\lesssim\;\rho^{\theta_0\tau}\delta(2\rho)+\Vert\nabla\bn_g\Vert_{L^2(D_{2\rho}(0)\setminus D_{\rho/2}(0))}\:.
\ees
It is not difficult to see that such an inequality implies that $\lim_{\rho\searrow0}\delta(\rho)=0$. Injecting this information back into the latter now gives
\bes
\rho^{-1}\delta^2(\rho)\;\lesssim\;\rho^{2\theta_0\tau-1}+\rho^{-1}\Vert\nabla\bn_g\Vert^2_{L^2(D_{2\rho}(0)\setminus D_{\rho/2}(0))}\:,
\ees
whence, using Fubini's theorem,
\bes
\int_0^{1/2}\rho^{-1}\delta^2(\rho)\,d\rho\;\lesssim\;1+\Vert\nabla\bn_g\Vert^2_{L^2(D_{1}(0))}\;<\;\infty\:,
\ees
as announced.\\[-3ex]

$\hfill\blacksquare$ \\

Since $|y|:=|\bp|(x)\simeq|x|^{\theta_0}$, the latter is equivalent to the desired
\bes
|\bA^g_{\bp}|(y)\;\lesssim\;|y|^{-1}\mu(|y|)\qquad\forall\:\:y\in\bp(\di)\quad\text{with}\quad|y|\ll1\:,
\ees
with $\lim_{|y|\searrow0}\mu(|y|)=\lim_{|x|\searrow0}\delta(|x|^{\theta_0})=0$.
This concludes the proof of Theorem \ref{eps}.

\subsubsection{The asymptotically synchronized case and the proof of Theorem \ref{synchroth}}

Throughout this section, we will suppose that $\tau$ is related to the integer $\theta_0$ in such a way that $\tau>1-1/{\theta_0}$. \\

Let us rewrite (\ref{equa2}) in the equivalent form
\bes
\text{div}\big(-\nabla\bH^g_{\bp}+2\pi_{T_g}\nabla\bH^g_{\bp}+|\bH^g_{\bp}|^2\nabla\bp+\bu_1)\;=\;\bu_2\qquad\text{on}\:\:\di\setminus\{0\}\:.
\ees
Using Lemma \ref{epsss}, it is not difficult to verify that (\ref{hypo27}) gives
\be\label{hypo2}
\left\{\begin{array}{lcl}
\bu_1&=&\text{O}\big(|x|^{\theta_0(\tau-1)-1}\delta(|x|)\big) \\[1.5ex]
\bu_2&=&\text{O}\big(|x|^{\theta_0(\tau-1)-2}\delta(|x|)\big)\:.
\end{array}\right.
\ee
On the other hand, we have also from Lemma \ref{epsss}:
\begin{eqnarray}\label{achier}
\pi_{T_g}\partial_{x^j}\bH^g_{\bp}&=&\pi_{T_g}\D{g}{\partial_{x^j}\bp}{}\bH^g_{\bp}\,+\text{O}\big(|\bp|^{-1+\tau}|\nabla\bp||\vec{H}^g_{\bp}||\big)
\nonumber\\[1ex]
&=&-\,\sum_{k=1,2}\big\langle\bH^g_{\bp}\,,\big(\vec{A}^g_{\bp}\big)_{j}^{k}\big\rangle_g\,\partial_{k}\bp\,+\,\text{O}\big(|x|^{\theta_0(\tau-1)-1}\delta(|x|)\big)\nonumber\\[1ex]
&=&\text{O}\big(|x|^{-\theta_0-1}\delta(|x|)\big)\:.
\end{eqnarray}
Owing to the latter and to (\ref{hypo2}), we may thus recast (\ref{equa2}) in the form
\bes
\text{div}\big(\nabla\bH^g_{\bp}+\vec{v}_1\big)\;=\;-\,\vec{u}_2\qquad\text{on}\:\:\:D_1(0)\setminus\{0\}\:,
\ees
where
\bes
\vec{v}_1\;:=\;-\,2\pi_{T_g}\nabla\bH^g_{\bp}-|\bH^g_{\bp}|^2\nabla\bp+\vec{u}_1\;=\;\text{O}\big(|x|^{-\theta_0-1}\delta(|x|)\big)\:.
\ees
As seen in Lemma \ref{epsss}, $|x|^{-1}\delta(|x|)$ is square integrable. It then follows that $|x|^{\theta_0}\vec{v}_1$ lies in $L^2(\di)$. \\
As for $\vec{u}_2$, it is such that $|x|^{1+\theta_0(1-\tau)}\bu_2$ lies in $L^2(\di)$. For notational convenience, we switch to the complex notation and replace the coordinates $(x^1,x^2)$ by the complex number $z$, in the usual way. Note that for some positive $\eta_1$ and $\eta_2$, we have
\bes
\big|z^{(1+\theta_0(1-\tau))/2}\vec{u}_2\big|\;\equiv\;|z|^{-(1+\theta_0(1-\tau))/2}|z|^{1+\theta_0(1-\tau)}|\vec{u}_2|\;\in\;L^{2+\eta_1}\cdot L^{2}\;\subset\;L^{1+\eta_2}\:,
\ees
where we have used the synchronization hypothesis $\tau>1-1/\theta_0$. We may thus introduce a Hodge decomposition 
\bes
\partial_{\bar{z}}\vec{w}_2\;=\;z^{(1+\theta_0(1-\tau))/2}\vec{u}_2\qquad\text{on}\:\:D_1(0)
\ees
and find that $\vec{w}_2$ lies in $L^{2+\eta_3}$, for some $\eta_3>0$. Hence, we have
\bes
\big|z^{-(1+\theta_0(1-\tau))/2}\vec{w}_2\big|\;\equiv\;|z|^{-(1+\theta_0(1-\tau))/2}|\vec{w}_2|\;\in\;L^{2+\eta_1}\cdot L^{2+\eta_3}\;\subset\;L^{1+\eta_4}\:,
\ees
for some $\eta_4>0$. We again perform a Hodge decomposition
\bes
\partial_{{z}}\vec{v}_2\;=\;-\,z^{-(1+\theta_0(1-\tau))/2}\vec{w}_2\qquad\text{on}\:\:D_1(0)
\ees
and find that the (necessarily real-valued) $\vec{v}_2$ satisfies
\bes
-\,\Delta\vec{v}_2\;=\;\vec{u}_2\qquad\text{on}\:\:D_1(0)\setminus\{0\}\:.
\ees
Moreover, since $\theta_0\ge1$, we have that
\be\label{chienquibraille}
|x|^{\theta_0}|\nabla\vec{v}_2|\;\equiv\;|z|^{\theta_0}|\partial_z\vec{v}_2|\;=\;|z|^{(-1+\theta_0(1+\tau))/2}|\vec{w}_2|\;\in\;L^\infty\cdot L^{2+\eta_3}\;\subset\;L^{2+\eta_3}\:,
\ee
where we have used that
\bes
-1+\theta_0(1+\tau)\;>\;-1+\theta_0(2-1/\theta_0)\;=\;2(\theta_0-1)\;\ge\;0\:,
\ees
which follows again from the synchronization hypothesis. \\[1.5ex]
Altogether, using the fact that $\theta_0\ge1$, the function $\bH^g_{\bp}$ satisfies a problem of the type
\bes
\text{div}\big(\nabla\bH^g_{\bp}+\bV\big)\;=\;\vec{0}\qquad\text{on}\:\:\:D_1(0)\setminus\{0\}\:,
\ees
where $\bV:=\vec{v}_1-\nabla\vec{v}_2$ satisfies $|x|^{\theta_0}\bV\in L^2$. In addition, we know that $|x|^{\theta_0-1}|\bH^g_{\bp}|$ lies as well in $L^2$.
According to Proposition \ref{appdev} in the appendix, we deduce that
\be\label{gradH}
|x|^{\theta_0}\nabla\bH^g_{\bp}\:\in\,L^2(D_1(0))\:.
\ee
For the record, let us note that (\ref{chienquibraille}) gives that $|x|^{\theta_0-1}\nabla\vec{v}_2$ lies in $L^{1+\eta_0}$ for some $\eta_0>0$ chosen small enough. \\

For the sake of our future needs, it is necessary to recast (\ref{equa2}) once more in a slightly more manageable form, namely
\bes
\text{div}\big(\nabla\bH^g_{\bp}-2\pi_{\bn_0}\nabla\vec{H}^g_{\bp}+|\bH^g_{\bp}|^2\nabla\bp+\vec{u}\big)\;=\;\vec{0}\qquad\text{on}\:\:\:D_1(0)\setminus\{0\}\:,
\ees
where
\bes
\vec{u}\;:=\;\vec{u}_1-\nabla\vec{v}_2+2(\pi_{\bn_0}-\pi_{\bn_g})\nabla\vec{H}^g_{\bp}\:.
\ees
We have just seen that $|x|^{\theta_0-1}\nabla\vec{v}_2$ lies in $L^{1+\eta_0}$ for some $\eta_0>0$. Furthermore, from our previous computations and (\ref{projn}), we find that
\begin{eqnarray}\label{defu}
|x|^{\theta_0-1}\big|\vec{u}+\nabla\vec{v}_2\big|&\lesssim&|x|^{\theta_0\tau-2}\delta(|x|)+|x|^{\theta_0(\tau+1)-1}|\nabla\bH^g_{\bp}|\nonumber\\[1ex]
&\lesssim&|x|^{\theta_0\tau-1}\big(|x|^{-1}\delta(|x|)+|x|^{\theta_0}|\nabla\bH^g_{\bp}|\big)\:,
\end{eqnarray}
which, we have shown, lies in the product of $L^{2+\eta}$ and of $L^2$, for some $\eta>0$. It then follows that
\be\label{hyper1}
|x|^{\theta_0-1}\vec{u}\:\in\;L^{1+\eta_0}(D_1(0))\qquad\text{for some}\:\:\eta_0>0\:.
\ee
An analogous argument reveals that
\be\label{hyper2}
|x|^{\theta_0}\vec{u}\:\in\;L^{2}(D_1(0))\:.
\ee

\medskip

We will now proceed studying (\ref{equa2}) in further details. To do so, we begin by defining the following constant vector called {\it residue}:
\bes
\vec{\gamma}_0\;:=\;\int_{\partial D_1(0)}\vec{\nu}\cdot\big(\nabla\bH^g_{\bp}-2\pi_{\bn_0}\nabla\vec{H}^g_{\bp}+|\bH^g_{\bp}|^2\nabla\bp+\vec{u}\big)
\ees
where $\vec{\nu}$ is the outward unit-normal to the flat unit-disk $D_1(0)$, and the dot product is understood, as always, as the standard Euclidean product in $\R^m$. \\
The equation (\ref{equa2}) implies that for any disk $D_\rho(0)$ of radius $\rho$ centered on the origin and contained in $D_1(0)\setminus\{0\}$, there holds
\bes
\int_{\partial D_\rho(0)} \vec{\nu}\cdot\big(\nabla\bH^g_{\bp}-2\pi_{\bn_0}\nabla\vec{H}^g_{\bp}+|\bH^g_{\bp}|^2\nabla\bp+\vec{u}\big)
\:=\:4\pi\,\vec{\gamma}_0\qquad\forall\:\:\rho\in(0,1)\:.
\ees
An elementary computation shows that
\bes
\int_{\partial D_{\rho}(0)} \vec{\nu}\cdot\nabla\log|x|\:=\:2\pi\:,\qquad\forall\:\:\rho>0\:.
\ees
Thus, upon setting
\bes
\bX\;:=\;\nabla\bH^g_{\bp}-2\pi_{\bn_0}\nabla\vec{H}^g_{\bp}+|\bH^g_{\bp}|^2\nabla\bp+\vec{u}\,-\,2\,\vec{\gamma}_0\,\nabla\log|x|\:,
\ees
we find
\bes
\text{div}\,\bX=0\quad\:\:\text{on}\:\:D_1(0)\setminus\{0\}\:,\hspace{.7cm}\text{and}\hspace{.5cm}\int_{\partial D_\rho(0)}\vec{\nu}\cdot\bX\,=\,{0}\qquad\forall\:\:\rho\in(0,1)\:.
\ees
As $\bX$ is smooth away from the origin, the Poincar\'e lemma implies the existence of an element $\bL\in C^{\infty}(D_1(0)\setminus\{0\})$, defined up to an additive constant, such that
\bes
\bX\;=\;\nabla^\perp\bL\;:=\;(-\,\partial_{x^2}\bL\,,\partial_{x^1}\bL)\qquad\text{on}\:\:\:\Om\:.
\ees
Note that Lemma \ref{epsss} yields
\bes
|x|^{\theta_0}|\bH^g_{\bp}|^2|\nabla\bp|(x)\;\lesssim\;|x|^{-1}\delta^2(|x|)\:\in\,L^2(D_1(0))\:.
\ees
From this, (\ref{gradH}), and (\ref{hyper2}), we deduce that $|x|^{\theta_0}\nabla\bL$ belongs to $L^2(D_1(0))$. 
A classical Hardy-Sobolev inequality gives the estimate
\be\label{estimL}
\theta_0^2\int_{D_1(0)}|x|^{2(\theta_0-1)}|\bL|^2\,dx\:\le\:\int_{D_1(0)}|x|^{2\theta_0}|\nabla\bL|^2\,dx\;+\;\theta_0\!\int_{\partial D_1(0)}|\bL|^2\;<\;\infty\:. 
\ee
The immersion $\bP$ has near the origin the asymptotic behavior $\,|\nabla\bp(x)|\simeq|x|^{\theta_0-1}$. Hence (\ref{estimL}) yields that 
\be\label{xx2}
\bL\cdot\nabla\bP\:,\;\bL\wedge\nabla\bP\;\in\;L^2(D_1(0))\:.
\ee

Next, we compute
\begin{eqnarray}\label{proS}
-\,\text{div}(\bL\cdot\nabla^\perp\bp)&=&\nabla\bp\cdot\nabla^\perp\bL\nonumber\\[1ex]
&=&\nabla\bp\cdot\nabla\bH^g_{\bp}+|\bH^g_{\bp}|^2|\nabla\bp|^2+\big(\vec{u}-2\,\vec{\gamma}_0\,\nabla\log|x|\big)\cdot\nabla\bp\nonumber\\[1ex]
&=&\text{div}(\bH^g_{\bp}\cdot\nabla\bp)+\big(|\nabla\bp|^2-|\nabla\bp|^2_g\big)|\bH^2_{\bp}|^2+\big(\vec{u}-2\,\vec{\gamma}_0\,\nabla\log|x|\big)\cdot\nabla\bp\:,
\end{eqnarray}
where we have used (\ref{laseulefaute}). \\
Let $f$ be the solution of 
\be\label{deff}
\left\{\begin{array}{rcllll}
\Delta f&=&\big(|\nabla\bp|^2-|\nabla\bp|^2_g\big)|\bH^2_{\bp}|^2+\big(\vec{u}-2\,\vec{\gamma}_0\,\nabla\log|x|\big)\cdot\nabla\bp&\quad&\text{in}\:\:\:D_1(0)\\[1ex]
f&=&0&\quad&\text{on}\:\:\:\partial D_1(0)\:.
\end{array}\right.
\ee
According to the asymptotic behavior of the metric near the origin, to Lemma \ref{epsss}, and to (\ref{hyper1}), we have
\bes
|\Delta f|\;\lesssim\;|x|^{\theta_0\tau-2}\delta(|x|)+|x|^{\theta_0-1}|\vec{u}|+|x|^{\theta_0-2}\:\in\,L^{1+\eta_0}(D_1(0))\quad\text{for some}\:\:\eta_0>0\:,
\ees
so that, in particular,
\be\label{regf}
\nabla f\:\in\,L^{2+\eta}(D_1(0))\qquad\text{for some}\:\:\eta>0\:.
\ee
For our future needs, we note that (\ref{proS}) states
\be\label{proS2}
\text{div}\big(\bL\cdot\nabla^\perp\bp+\bH^g_{\bp}\cdot\nabla\bp+\nabla f  \big)\;=\;0\qquad\text{in}\:\:\Om\:.
\ee
Similarly, again using (\ref{laseulefaute}), we now compute
\begin{eqnarray}\label{proR}
-\,\text{div}(\bL\wedge\nabla^\perp\bp)&=&\nabla\bp\wedge\nabla^\perp\bL\nonumber\\[1ex]
&=&\nabla\bp\wedge\nabla\bH^g_{\bp}-2\nabla\bp\wedge\pi_{\bn_0}\nabla\bH^g_{\bp}-\big(\vec{u}-2\,\vec{\gamma}_0\,\nabla\log|x|\big)\wedge\nabla\bp\nonumber\\[1ex]
&=&\text{div}(\bH^g_{\bp}\wedge\nabla\bp)+\bF_1+2\nabla\bp\wedge\pi_{T_0}\nabla\bH^g_{\bp}-\big(\vec{u}-2\,\vec{\gamma}_0\,\nabla\log|x|\big)\wedge\nabla\bp\:,
\end{eqnarray}
where it is easy to check from (\ref{laseulefaute}) that for some $\eta>0$:
\be\label{chiantox}
|\bF_1|\;=\;\text{O}\big(|\bp|^{\tau-1}|\nabla\bp|^2||\vec{H}^g_{\bp}|  \big)\;=\;|x|^{\theta_0\tau-1}\text{O}\big(\nabla\bp||\vec{H}^g_{\bp}|  \big)\,\in\,L^{1+\eta}\:.
\ee
This will be used shortly.\\
Previously encountered estimates give
\begin{eqnarray}\label{ctemerde}
\nabla\bp\wedge\pi_{T_0}\nabla\bH^g_{\bp}&=&\nabla\bp\wedge\pi_{T_g}\nabla\bH^g_{\bp}+\text{O}\big(|\nabla\bp||\nabla\bH^g_{\bp}||\bp|^{\tau}\big)\nonumber\\[1ex]
&=&\partial_{x^1}\bp\wedge\pi_{T_g}\nabli{g}{\partial_{x^1}\bp}\bH^g_{\bp}+\partial_{x^2}\bp\wedge\pi_{T_g}\nabli{g}{\partial_{x^2}\bp}\bH^g_{\bp}\nonumber\\
&&\hspace{1.5cm}+\:\text{O}\big(|\nabla\bp||\nabla\bH^g_{\bp}||\bp|^{\tau}+|\nabla\bp|^2|\bH^g_{\bp}||\bp|^{\tau-1}\big)\nonumber\\[1ex]
&=&|x|^{\theta_0\tau-1}\,\text{O}\big(|x|^{\theta_0}|\nabla\bH^g_{\bp}|+|x|^{\theta_0-1}|\bH^g_{\bp}|\big)\:,
\end{eqnarray}
where we have used the easily-verified fact that
\bes
\partial_{x^1}\bp\wedge\pi_{T_g}\nabli{g}{\partial_{x^1}\bp}\bH^g_{\bp}+\partial_{x^2}\bp\wedge\pi_{T_g}\nabli{g}{\partial_{x^2}\bp}\bH^g_{\bp}\;=\;\vec{0}
\ees
which follows from the symmetry of the second-fundamental form.\\
According to (\ref{gradH}), the bracketed term on the right-hand side of (\ref{ctemerde}) lies in $L^2$. In addition, the factor $|x|^{\theta_0\tau-1}$ surely lies in $L^{2+\eta'}$, for some suitably chosen $\eta'>0$. It then follows that
\be\label{hyper3}
\nabla\bp\wedge\pi_{T_0}\nabla\bH^g_{\bp}\:\in\;L^{1+\eta_0}(D_1(0))\qquad\text{for some}\:\:\eta_0>0\:.
\ee
Let now $\bF$ be the solution of 
\bes
\left\{\begin{array}{rcllll}
\Delta\bF&=&2\nabla\bp\wedge\pi_{T_0}\nabla\bH^g_{\bp}-\big(\vec{u}-2\,\vec{\gamma}_0\,\nabla\log|x|\big)\wedge\nabla\bp+\bF_1&\quad&\text{in}\:\:\:D_1(0)\\[1ex]
\bF&=&\vec{0}&\quad&\text{on}\:\:\:\partial D_1(0)\:.
\end{array}\right.
\ees
With the help of (\ref{hyper1}), (\ref{chiantox}), and (\ref{hyper3}), we have that $\Delta\bF$ lies in $L^{1+\eta_0}(\di)$ for some $\eta_0>0$. Hence,
\be\label{regF}
\nabla\bF\:\in\,L^{2+\eta}(D_1(0))\qquad\text{for some}\:\:\eta>0\:.
\ee
For our future needs, we note that (\ref{proR}) states
\be\label{proR2}
\text{div}\big(\bL\wedge\nabla^\perp\bp+\bH^g_{\bp}\wedge\nabla\bp+\nabla\bF \big)\;=\;\vec{0}\qquad\text{in}\:\:\Om\:.
\ee

Note that the terms under the divergence symbols in (\ref{proS2}) and in (\ref{proR2}) both belong to $L^2(D_1(0))$, owing to (\ref{gradH}), (\ref{xx2}), (\ref{regf}), and to (\ref{regF}). The distributional equations (\ref{proS2}) and (\ref{proR2}), which are {\it a priori} to be understood on $\Om$, may thus be extended to all of $D_1(0)$. Indeed, a classical result of Laurent Schwartz states that the only distributions supported on $\{0\}$ are linear combinations of derivatives of the Dirac delta mass. Yet, none of these (including delta itself) belongs to $W^{-1,2}$. We shall thus understand (\ref{proS2}) and (\ref{proR2}) on $D_1(0)$. 
It is not difficult to verify (cf. Corollary IX.5 in \cite{DL}) that a divergence-free vector field in $L^2(D_1(0))$ is the curl of an element in $W^{1,2}(D_1(0))$. 
We apply this observation to (\ref{proS2}) and in (\ref{proR2}) so as to infer the existence of two functions\footnote{$S$ is a scalar while $\bR$ is $\bigwedge^2(\R^m)$-valued.} $S$ and of $\bR$ in the space $W^{1,2}(D_1(0))\cap C^\infty(\Om)$, with
\bes
\left\{\begin{array}{rclll}
\nabla^\perp S&=&\bL\cdot\nabla^\perp\bp+\bH^g_{\bp}\cdot\nabla\bp+\nabla f&&\\[1.5ex]
\nabla^\perp\bR&=&\bL\wedge\nabla^\perp\bp+\bH^g_{\bp}\wedge\nabla\bp+\nabla\bF \:.
\end{array}\right.
\ees
According to Lemma \ref{identities} from the Appendix, the functions $S$ and $\bR$ satisfy on $D_1(0)$ the following equations:
\be\label{sysSR-}
\left\{\begin{array}{rclll}
-\,\nabla S&=&\nabla^\perp f + (\star_g\bn_g)\cdot(\nabla^\perp\bR-\nabla\bF)+q &&\\[1.5ex]
-\,\nabla\bR&=&\nabla^\perp\bF+(\star_g\bn_g)\bul(\nabla^\perp\bR-\nabla\bF)-(\star_g\bn_g)(\nabla^\perp S-\nabla f)+\bQ\:,
\end{array}\right.
\ee
where 
\be\label{hypq}
|q|+|\bQ|\;=\;\text{e}^{\la}\big(|\bL|+|\bH^g_{\bp}|\big)\text{O}_2\big(|\bp|^{\tau}\big)\;=\;\text{O}\big(|x|^{\theta_0(\tau+1)-1}\big)\big(|\bL|+|\bH^g_{\bp}|\big)\:.
\ee
Note that 
\bes
\big|\nabla\big(|x|^{\theta_0(\tau+1)-1}\bL\big)\big|\;\leq\;|x|^{\theta_0\tau-1}\big(|x|^{\theta_0-1}|\bL|+|x|^{\theta_0}|\nabla\bL|\big)\:.
\ees
As we have already oftentimes seen, the first factor on the right-hand side lies in $L^{2+\eta'}$, for some $\eta'>0$, while the second factor on the right-hand side of the latter belongs to $L^2$. Accordingly, $|x|^{\theta_0(\tau+1)-1}\bL\in W^{1,1+\eta_0}$ for some $\eta_0>0$, from which it follows that $|x|^{\theta_0(\tau+1)-1}\bL\in L^{2+\eta}$ for some $\eta>0$. For the exact same reason, we have $|x|^{\theta_0(\tau+1)-1}\bH^g_{\bp}\in L^{2+\eta}$ for some $\eta>0$. Bringing this into (\ref{hypq}) shows that
\be\label{hypq2}
|q|+|\bQ|\;\in\,L^{2+\eta}\:.
\ee 
Differentiating (\ref{sysSR-}) throughout yields
\be\label{sysSR+}
\left\{\begin{array}{rclll}
-\,\Delta S&=&\nabla(\star_g\bn_g)\cdot\nabla^\perp\bR\,-\,\text{div}\big((\star_g\bn_g)\cdot\nabla\bF+q \big) &&\\[1.5ex]
-\,\Delta\bR&=&\nabla(\star_g\bn_g)\bul\nabla^\perp\bR\,-\,\nabla(\star_g\bn_g)\cdot\nabla^\perp S\\[1ex]
&&\hspace{2cm}\,-\,\text{div}\big((\star_g\bn_g)\bul\nabla\bF-(\star_g\bn_g)\nabla f+\bQ \big)\:.
\end{array}\right.
\ee
From (\ref{regf}), (\ref{regF}), and (\ref{hypq2}), the terms under the divergence forms on the right-hand side belong to $L^{2+\eta}$ for some $\eta>0$. On the other hand, we have seen that $\nabla S$ and $\nabla\bR$ lie in $L^2$. And finally, (\ref{epshyp}) guarantees that the $L^2$-norm of $\nabla(\star_g\bn_g)$ may be chosen as small as we please. We are thus in the position of applying Proposition \ref{adams} from the Appendix to conclude that there exists $p>2$ such that
\be\label{betterSR}
\nabla S\:,\:\nabla\bR\;\in\,L^p(\di)\:.
\ee

We learn in Lemma \ref{identities} that
\be\label{oubli}
2\text{e}^{2\la}\bH^g_{\bp}\;=\;(\nabla S+\nabla^\perp f)\cdot\nabla\bp+(\nabla\bR+\nabla^\perp\bF)\bul\nabla^\perp\bp+\text{e}^{2\la}|\bH^g_{\bp}|\text{O}_2(|\bp|^{\tau})\:.
\ee
Using the known asymptotic behaviors of $\bp$ and of its gradient near the origin, along with (\ref{laseulefaute}),  the latter reads
\bes
\Delta\bp\;=\;\Big[(\nabla S+\nabla^\perp f)\cdot\nabla\bp+(\nabla\bR+\nabla^\perp\bF)\bul\nabla^\perp\bp\Big]\,\big(1+\text{O}(|x|^{\theta_0\tau})\big)+\text{O}\big(|x|^{\theta_0(\tau+1)-2}\big)\:,
\ees
so that
\be\label{tuz}
\text{e}^{-\la}|\Delta\bp|\;\leq\;\Big(|\nabla S|+|\nabla\bR|+|\nabla f|+|\nabla\bF|\Big)\big(1+\text{O}(|x|^{\theta_0\tau})\big)+\text{O}\big(|x|^{\theta_0\tau-1}\big)\:,
\ee
where we have used Lemma \ref{epsss}.\\
Owing to (\ref{regf}), (\ref{regF}), and to (\ref{betterSR}), we see that the right-hand side of the latter belongs to $L^{p}(\di)$ for some $p>2$. We may thus call upon Proposition \ref{CZ} from the Appendix to conclude that near the origin, the immersion $\bp$ displays an asymptotic behavior of the form:
\bes
(\partial_{x^1}\!+i\,\partial_{x^2})\bp(x)\;=\;\bPe(\overline{x})+\,|x|^{\theta_0-1}\bT(x)\:,
\ees
where $\bPe$ is a $\C^m$-valued polynomial of degree at most $(\theta_0-1)$, and $\vec{T}(x)=\text{O}\big(|x|^{1-\frac{2}{p}-\epsilon'}\big)$ for every $\epsilon'>0$. Because $\,\text{e}^{-\la}\nabla\bp\,$ is a bounded function, we deduce more precisely that $\bPe(\overline{x})=\theta_0\bB^*\overline{x}^{\,\theta_0-1}$, for some constant vector $\bB\in\C^m$\, (we denote its complex conjugate by $\bB^*$), so that
\bes
\nabla\bp(x)\;=\;\left(\begin{array}{c}\Re\\[.5ex]-\,\Im\end{array}\right)\big(\theta_0\vec{B}{x}^{\theta_0-1}\big)\,+\,|x|^{\theta_0-1}\bT(x)\:.
\ees
Equivalently, switching to the complex notation, there holds
\be\label{locexphi}
\pz\bp\;=\;\dfrac{\theta_0}{2}\bB z^{\theta_0-1}+\,\text{O}\big(|z|^{\theta_0-\frac{2}{p}-\epsilon'}\big)\quad\forall\:\epsilon'>0\:.
\ee
We write $\bB=\bB_R+i\bB_I\in\R^2\otimes\R^m$. 
The conformality condition on $\bp$ shows easily that $|\bB|^2_g=0$, whence
\be\label{stuffonb}
|\bB_{R}|_g\;=\;|\bB_{I}|_g\qquad\text{and}\qquad\langle\bB_{R},\bB_{I}\rangle_g\;=\;0\:.
\ee
Yet more precisely, as $|\nabla\bp|_g^2=2\,\text{e}^{2\la}$, we see that
\bes
|\bB_{R}|_g\;=\;|\bB_{I}|_g\;=\;\dfrac{1}{\theta_0}\,\lim_{z\rightarrow0}\,\dfrac{\text{e}^{\la(z,\bar{z})}}{|z|^{\theta_0-1}}\,\in\:\:]0\,,\infty[\:.
\ees
Because $\bp(0)=\vec{0}$, we obtain from (\ref{locexphi}) the local expansion
\be\label{exphi}
\bp(z,\bar{z})\;=\;\Re\big(\bB z^{\theta_0}\big)\,+\,\text{O}_1\big(|z|^{\theta_0+1-\frac{2}{p}-\epsilon'}\big)\:.
\ee
On the other hand, from $\pi_{\bn_g}\!\nabla\bp\equiv\vec{0}$, we deduce from (\ref{locexphi}) that
\bes
\pi_{\bn_g}\bB\;=\;\text{O}\big(|z|^{1-\frac{2}{p}-\epsilon'}\big)\qquad\:\:\:\forall\:\epsilon'>0\:.
\ees
Let now $\delta:=1-\frac{2}{p}\in(0,1)$, and let $0<\eta<p$ be arbitrary. We choose some $\epsilon'$ satisfying
\bes
0\;<\;\epsilon'\;<\;\frac{2\,\eta}{p(p-\eta)}\;\equiv\;\delta-1+\frac{2}{p-\eta}\:.
\ees
We have observed that $\,\pi_{\bn_g}\bB=\text{O}(|z|^{\delta-\epsilon'})$, hence $\,
\pi_{\bn_g}\bB=\text{o}\big(|z|^{1-\frac{2}{p-\eta}}\big)$\,,
and in particular, we find
\be\label{pina}
|z|^{-1}\pi_{\bn_g}\bB\;\in\;L^{p-\eta}(\di)\:\:\:\:\qquad\forall\:\:\eta>0\:.
\ee
This fact shall be put to good use in the sequel.\\

Proposition \ref{introms} states that the weight $\text{e}^{\la}$ satisfies the conditions of  Proposition \ref{CZ}-(ii). Hence, we deduce from (\ref{tuz}) that
\be\label{locexdelphi}
\nabla^2\bp\;=\;\theta_0\,(1-\theta_0)\left(\begin{array}{cc}-\,\Re&\Im\\[.5ex]\Im&\Re\end{array}\right)\big(\vec{B}{z}^{\theta_0-2}\big)\,+\,|x|^{\theta_0-1}\vec{Z}\:,
\ee
where $\vec{B}$ is as in (\ref{locexphi}), and $\vec{Z}$ lies in $\mathbb{R}^4\otimes L^{p-\epsilon'}(\di,\R^m)$ for every $\epsilon'>0$. The exponent $p>2$ is the same as above. We obtain from (\ref{locexdelphi}) that
\bes
\text{e}^{-\la}\big|\pi_{\bn_g}\!\nabla^2\bP\big|\;\lesssim\;|z|^{-1}|\pi_{\bn_g}\vec{B}|\,+\,|\pi_{\bn_g}\vec{Z}|\:.
\ees
According to (\ref{pina}), the first summand on the right-hand side of the latter belongs to $L^{p-\eta}\,$ for all $\eta>0$. Moreover, we have seen that $\pi_{\bn_g}\vec{Z}$ lies in $L^{p-\epsilon'}$ for all $\epsilon'>0$. Whence, it follows that $\,\text{e}^{-\la}\pi_{\bn_g}\!\nabla^2\bP\,$ is itself an element of $L^{p-\epsilon'}$ for all $\epsilon'>0$. By definition, this confirms that the regularity of the second fundamental form has been improved to
\be\label{regn}
\text{e}^{\la}\!\vec{A}^g_{\bp}\;\in\;L^{p-\epsilon'}(\di)\:,\qquad\quad\forall\:\:\epsilon'>0\:.
\ee
As we have seen, $\nabla(\star_g\bn_g)$ inherits the integrability of $\text{e}^{\la}\!\vec{A}^g_{\bp}$, so that
\be\label{regnn}
\nabla(\star_g\bn_g)\;\in\;L^{p-\epsilon'}(\di)\:,\qquad\quad\forall\:\:\epsilon'>0\:.
\ee
In the sequel, we will fix $t:=p-\epsilon'>2$. 
In light of this new fact, we return to the proof of Lemma \ref{epsss} to find that the function $\delta(r)$ defined in (\ref{deldef}) now satisfies
\be\label{celleci}
\delta(r)\;\lesssim\;r^{1-\frac{2}{t}}\qquad\text{and}\qquad |x|^{-1}\delta(|x|)\in\,L^{t}(\di)\:.
\ee
Having this information at our disposal, it is not difficult to follow the stream of our previous argument and to successively find that 
\bes
|x|^{\theta_0-1}\nabla\vec{v}_2\;\;,\;\;|x|^{\theta_0}\nabla\bH^g_{\bp}\:\:,\:\:|x|^{\theta_0}|\bH^g_{\bp}|^2|\nabla\bp|\:\:,\:\: |x|^{\theta_0}\nabla\bL\:\:,\:\: |x|^{\theta_0-1}\bL\;\in\,L^t(\di)\:.
\ees 
From this and an argument analogous to (\ref{defu}), we also obtain that
\be\label{regu2}
|x|^{\theta_0-1}\vec{u}\;\in\;\left\{\begin{array}{rcl}L^{s}(\di)&,&\theta_0=1\\[1ex]
L^t(\di)&,&\theta_0\ge2\:,  
\end{array}\right.
\ee
where, for our convenience, we choose $s>1$ to be such that\footnote{this requires that $\epsilon'$ be chosen small enough.}
\bes
\dfrac{1}{s}\;=\;\dfrac{1}{t}+\dfrac{1}{\sigma}\qquad\text{for any}\qquad0\,<\,\dfrac{1}{\sigma}\,<\,\dfrac{1}{p}-\dfrac{1}{t}+\dfrac{1}{2}\:.
\ees
Introduced into (\ref{deff}) this yields in turn that\footnote{The weak-$L^2$ Marcinkiewicz space $L^{2,\infty}(\di)$ is defined as those functions $f$ with the property that $\:\sup_{\alpha>0}\alpha^2\Big|\big\{x\in \di\,;\,|f(x)|\ge\alpha\big\}\Big|<\infty$. In dimension two, the prototype element of $L^{2,\infty}$ is $|x|^{-1}\,$.}
 \be\label{regf2}
\nabla^2 f\;\in\;\left\{\begin{array}{lcl}L^{s}(\di)&,&\theta_0=1\:\:\:\text{and}\:\:\:s<2\\[1ex]
L^{2,\infty}(\di)&,&\theta_0=1\:\:\:\text{and}\:\:\:s\ge2\\[1ex]
L^t(\di)&,&\theta_0\ge2\:,
\end{array}\right.
\ee
whence\footnote{we also use a result of Luc Tartar \cite{Tar} stating that $\,W^{1,(2,\infty)}\subset{BMO}$.}
\bes
\nabla f\;\in\;\left\{\begin{array}{lcl}L^{s^*}&,&\theta_0=1\:\:\:\text{and}\:\:\:s<2\\[1ex]
BMO&,&\theta_0=1\:\:\:\text{and}\:\:\:s\ge2\\[1ex]
L^\infty&,&\theta_0\ge2\:.
\end{array}\right.
\ees
where
\bes
\dfrac{1}{s^*}\;:=\;\dfrac{1}{s}\,-\,\dfrac{1}{2}\;<\;\dfrac{1}{p}\:.
\ees
Differentiating throughout the order relation (\ref{hypq}) and using the fact that $|x|\nabla\la$ is bounded gives
\bes
|\nabla\bQ|+|\nabla q|\;\lesssim\;|x|^{\theta_0\tau-1}\Big[|x|^{\theta_0-1}(|\bL|+|\bH^g_{\bp}|)+|x|^{\theta_0}(|\nabla\bL|+|\nabla\bH^g_{\bp}|)  \Big]\:.
\ees
Recall that we are assuming $\theta_0\tau-1>\theta_0-2$ in this section. As the bracketed term on the right-hand side lies in $L^t$, we obtain that
\be\label{regq2}
\nabla q\:,\:\nabla\bQ\;\in\;\left\{\begin{array}{rcl}L^{s}(\di)&,&\theta_0=1\\[1ex]
L^t(\di)&,&\theta_0\ge2\:,  
\end{array}\right.
\ee
where the exponent $s$ is as in (\ref{regu2}). \\
Bringing (\ref{betterSR}), (\ref{regnn}), (\ref{regf2}), and (\ref{regq2}) into  (\ref{sysSR+}) shows that $\Delta S$ lies in one of three possible spaces, namely $L^p\cdot L^t$, $L^{2,\infty}$, $L^s$,
whichever one is the smallest. The Sobolev embedding theorem then gives that $\nabla S$ lies in the smallest among $L^r$, $BMO$, and $L^{s^*}$, where $s^*$ is as above, and
\bes
\dfrac{1}{r}\;=\;\dfrac{1}{p}+\dfrac{1}{t}-\dfrac{1}{2}\;=\;\dfrac{1}{p}+\dfrac{1}{p-\epsilon'}-\dfrac{1}{2}\;<\;\dfrac{1}{p}\:.
\ees
As we have seen above, $s^*>p$, and clearly $BMO$ contains all $L^p$ spaces for $p$ finite. Accordingly, in all configurations, we see that the integrability of $\nabla S$ has been improved. Identical reasoning and conclusion hold with $\bF$ and $\bR$ respectively in place of $f$ and $S$. This procedure may be repeated until reaching that 
\bes
\nabla S\:\,,\,\nabla\bR\:\,\in\:L^{b}(\di)\qquad\:\:\forall\:\:b<\infty\:.
\ees
Introducing this information back into the above procedure yields that
\bes
\nabla S\;\,,\,\nabla\bR\:\,\in\:\,\left\{\begin{array}{lcl}W^{1,(2,\infty)}(\di)&,&\text{if}\:\:\:\theta_0=1\\[1ex]
W^{1,b}(\di)&,&\text{if}\:\:\:\theta_0\ge2\:,\qquad\forall\:\:b<\infty\:.
\end{array}\right.
\ees

From (\ref{oubli}), we have that
\bes
|x|^{\theta_0-1}|\bH^g_{\bp}|\;\lesssim\;|\nabla S|+|\nabla f|+|\nabla\bR|+|\nabla\bF|\,\in\,BMO(\di)\:,
\ees
which, once fed back into (\ref{oubli}) gives that $|x|^{1-\theta_0}\Delta\bp\in\bigcap_{p<\infty}L^p$. We may thus once more call upon Proposition \ref{CZ} to obtain the following improvement of (\ref{exphi}):
\bes
\bp\;=\;\Re\big(\bB z^{\theta_0})+\text{O}_1\big(|z|^{\theta_0+1-\epsilon'}\big)\qquad\forall\:\epsilon'>0\:,
\ees
and, just as (\ref{regn}) and (\ref{celleci}) followed from (\ref{locexphi}), we get this time that 
\bes
\text{e}^{\la}|\vec{A}^g_{\bp}|\,\in\,\bigcap_{p<\infty}L^p\:,
\ees
and that the function $\delta(r)$ defined in (\ref{deldef}) now satisfies
\be\label{celleci2}
\delta(r)\;\lesssim\;r^{1-\epsilon'}\quad\forall\:\:\epsilon'>0\qquad\text{and}\qquad |x|^{-1}\delta(|x|)\in\,\bigcap_{t<\infty}L^{t}(\di)\:.
\ee

To finish the proof of Theorem \ref{synchroth}, we now need to distinguish two separate cases. 

\paragraph{\underline{The case $\theta_0\tau>1$.}} This requires in particular that $\theta_0\ge2$. The case $\theta_0=1$ will be handled separately. From (\ref{hypo27}) and (\ref{celleci2}), it easily follows that $|x|^{\theta_0}\bu_1$ and $|x|^{\theta_0}\bu_2$ lie in $L^p$ for all $p<\infty$. Proceeding exactly as before, we successively obtain that
\be\label{mosta}
|x|^{\theta_0-1}\nabla\vec{v}_2\;\;,\;\;|x|^{\theta_0}\nabla\bH^g_{\bp}\;\;,\;\;|x|^{\theta_0}|\bH^g_{\bp}|^2\nabla\bp\;\;,\;\;|x|^{\theta_0-1}\bL\:\in\,\bigcap_{p<\infty}L^p(\di)\:.
\ee

Let us return to the equation defining $\bL$, namely 
\be\label{defL2}
\nabla^\perp\bL\;=\;-\,\nabla\bH^g_{\bp}+2\pi_{T_g}\nabla\bH^g_{\bp}+|\bH^g_{\bp}|^2\nabla\bp+\vec{u}_1-\nabla\vec{v}_2-2\vec{\gamma}_0\nabla\log|x|\:.
\ee
For notational convenience, let us set
\bes
2\bJ\;:=\;2\pi_{T_g}\nabla\bH^g_{\bp}+|\bH^g_{\bp}|^2\nabla\bp+\vec{u}_1-\nabla\vec{v}_2\:,
\ees
and note that
\be\label{regJ}
|x|^{\theta_0-1}\vec{J}\,\in\,\bigcap_{p<\infty}L^{p}(\di)\:.
\ee
Introducing complex coordinates $z:=x^1+ix^2$ and $\bar{z}:=x^1-ix^2$, we may recast (\ref{defL2}) in the form
\bes
\pzb\big(i\bL+\bH^g_{\bp}+2\vec{\gamma}_0\log|z|\big)\;=\;\vec{J}\qquad\text{on}\:\:\:\Om\:.
\ees
Any function $\vec{w}$ satisfying 
\bes
\pzb\vec{w}\;=\;z^{\theta_0-1}\bJ\qquad\text{on}\:\:\:\di
\ees
lies in $C^{0,1-\epsilon'}(\di)$ for all $\epsilon'>0$, owing to (\ref{regJ}). Furthermore,
\be\label{tutu1}
\pzb\Big[z^{\theta_0-1}\big(i\bL+\bH^g_{\bp}+2\vec{\gamma}_0\log|z|\big)-(\vec{w}-\vec{w}(0))\Big]\;=\;\vec{0}\qquad\text{on}\:\:\:\Om\:.
\ee
From (\ref{mosta}), one sees that the bracketed function in the latter lies in $L^{2+\eta}(\di)$ for some (in fact for {\it all}) $\eta>0$. The equation (\ref{tutu1}) thus extends to all of the unit disk, and there exists some holomorphic function $\vec{E}$ such that
\bes
z^{\theta_0-1}\big(i\bL+\bH^g_{\bp}+2\vec{\gamma}_0\log|z|\big)-(\vec{w}-\vec{w}(0))\;=\;\vec{E}\:.
\ees
Hence, as $\vec{w}$ is H\"older continuous, 
\bes
i\bL+\bH^g_{\bp}\;=\;-\,2\vec{\gamma}_0\log|z|+\vec{E}_0\,z^{1-\theta_0}+\text{O}\big(|z|^{2-\theta_0-\epsilon'}\big)\qquad\forall\:\:\epsilon'>0\:,
\ees
for some constant $\bE_0\in\C^m$.
In particular, since $\bL$ is real-valued, we find
\bes
\bH^g_{\bp}(x)\;=\;-\,2\vec{\gamma}_0\log|x|+\Re\big(\vec{E}_0\,x^{1-\theta_0}\big)+\text{O}\big(|x|^{2-\theta_0-\epsilon'}\big)\qquad\forall\:\:\epsilon'>0\:.
\ees
Brought into the equation (\ref{laseulefaute}), namely $\Delta\bp=2\text{e}^{2\la}\bH^g_{\bp}+\text{O}(|x|^{\theta_0(\tau+1)-2})$, and using the fact that the weight $\text{e}^{\la}\simeq|x|^{\theta_0-1}$ satisfies the conditions of Propositions \ref{CZ}-(ii) (owing to (\ref{ms7})) yields now the local asymptotic expansion (valid for all $\epsilon'>0$):
\bes
\bp(x)\;=\;\Re\Big(\bB x^{\theta_0}+\bB_1x^{\theta_0+1}+\bB_2|x|^{2\theta_0}x^{1-\theta_0} \Big)+\text{O}_2\big(|x|^{\theta_0(\tau+1)-\epsilon'}+|x|^{\theta_0+2-\epsilon'}\big)\:.
\ees
Naturally, depending on the relative sizes of $\tau$ and $\theta_0$, one summand of the remainder term will be absorbed into the other.

\paragraph{\underline{The case $\theta_0\tau\leq1$.}} For an immersion synchronized with the ambient metric, we have $\tau>1-1/\theta_0$. Since $\theta_0\ge1$ is an integer, the only possibility is $\theta_0=1$. Without much difficulty, we verify, just as was previously done, that now $|x|\bu_1$ and $|x|\bu_2$ lie in $L^p(\di)$ for all $p<2/(1-\tau)$. Proceeding exactly as before, we successively obtain that
\be\label{omosta}
\nabla\vec{v}_2\;\;,\;\;|x|\nabla\bH^g_{\bp}\;\;,\;\;|x||\bH^g_{\bp}|^2\nabla\bp\;\;,\;\;\bL\:\in\,\bigcap_{p<\frac{2}{1-\tau}}L^p(\di)\:.
\ee
Let us return to the equation defining $\bL$, namely 
\be\label{odefL2}
\nabla^\perp\bL\;=\;-\,\nabla\bH^g_{\bp}+2\pi_{T_g}\nabla\bH^g_{\bp}+|\bH^g_{\bp}|^2\nabla\bp+\vec{u}_1-\nabla\vec{v}_2-2\vec{\gamma}_0\nabla\log|x|\:.
\ee
For notational convenience, let us set
\bes
2\bJ\;:=\;2\pi_{T_g}\nabla\bH^g_{\bp}+|\bH^g_{\bp}|^2\nabla\bp+\vec{u}_1-\nabla\vec{v}_2\:,
\ees
and note that
\be\label{oregJ}
\vec{J}\,\in\,\bigcap_{p<\frac{2}{1-\tau}}L^{p}(\di)\:.
\ee
Introducing complex coordinates $z:=x^1+ix^2$ and $\bar{z}:=x^1-ix^2$, we may recast (\ref{odefL2}) in the form
\bes
\pzb\big(i\bL+\bH^g_{\bp}+2\vec{\gamma}_0\log|z|\big)\;=\;\vec{J}\qquad\text{on}\:\:\:\Om\:.
\ees
Any function $\vec{w}$ satisfying 
\bes
\pzb\vec{w}\;=\;\bJ\qquad\text{on}\:\:\:\di
\ees
lies in $C^{0,\tau-\epsilon'}(\di)$ for all $\epsilon'>0$, owing to (\ref{oregJ}). Furthermore,
\be\label{otutu1}
\pzb\Big[i\bL+\bH^g_{\bp}+2\vec{\gamma}_0\log|z|-(\vec{w}-\vec{w}(0))\Big]\;=\;\vec{0}\qquad\text{on}\:\:\:\Om\:.
\ee
From (\ref{omosta}) and the fact that $\tau>0$, one sees that the bracketed function in the latter lies in $L^{2+\eta}(\di)$ for some $\eta>0$. The equation (\ref{otutu1}) thus extends to all of the unit disk, and there exists some holomorphic function $\vec{E}$ such that
\bes
i\bL+\bH^g_{\bp}+2\vec{\gamma}_0\log|z|-(\vec{w}-\vec{w}(0))\;=\;\vec{E}\:.
\ees
Hence, as $\vec{w}$ is H\"older continuous, 
\bes
i\bL+\bH^g_{\bp}\;=\;-\,2\vec{\gamma}_0\log|z|+\vec{E}_0\,z+\text{O}\big(|z|^{\tau-\epsilon'}\big)\qquad\forall\:\:\epsilon'>0\:,
\ees
for some constant $\bE_0\in\C^m$.
In particular, since $\bL$ is real-valued, we find
\be\label{zut2}
\bH^g_{\bp}(x)\;=\;-\,2\vec{\gamma}_0\log|x|+\Re\big(\vec{E}_0\big)+\text{O}\big(|x|^{\tau-\epsilon'}\big)\qquad\forall\:\:\epsilon'>0\:.
\ee
Brought into the equation (\ref{laseulefaute}), namely $\Delta\bp=2\text{e}^{2\la}\bH^g_{\bp}+\text{O}(|x|^{\tau-1})$, and using the fact that the weight $\text{e}^{\la}\simeq1$ yields now the local asymptotic expansion (valid for all $\epsilon'>0$):
\begin{equation}\label{serfs}
\bp(x)\;=\;\Re\big(\bB x\big)+\text{O}_2\big(|x|^{\tau+1-\epsilon'}\big)\:.
\end{equation}

\subsubsection{Embeddings in asymptotically Schwarzschild spaces and the proof of Theorem \ref{schwat}}\label{allld}

In this section, we will consider an ambient metric of Schwarzschild decay, namely
\be\label{schwatooo}
g_{\al\beta}(y)\;=\;\big(1+c|y|\big)\delta_{ij}+\text{O}_2(|y|^{1+\kappa})\:,\qquad|y|\ll1\:,
\ee
for some $\kappa>0$ (we will also without loss of generality suppose that $\kappa<1$) and some constant $c$. The Christoffel symbols of such a metric are easily computed. Moreover, we can use (\ref{serfs}) with $\tau=1$ to replace formula (\ref{laseulefaute}) by
\bes
\Delta\bp\;=\;c_1\bH^g_{\bp}+\vec{c}_2+\text{O}(|x|^{\kappa})\:,\qquad|x|\ll1\:,
\ees
for some nonzero constants $c_1\in\R$ and $\vec{c}_2\in\R^m$. 
Note that we have used that $\text{e}^{\la(x)}$ has a nonzero limit as $x$ approaches the origin. We have seen in the previous section that the mean curvature vector satisfies the local expansion (\ref{zut2}). Hence
\bes
\Delta\bp\;=\;-2\vec{\gamma}_0c_1\log|x|+\vec{E}_1+\text{O}(|x|^{\kappa})\:,\qquad|x|\ll1\:,
\ees
where we have set $\vec{E}_1:=c_1\Re(\vec{E}_0)+\vec{c}_2$. As we have previously done, we use Proposition \ref{CZ} to obtain the local expansion (valid for all $\epsilon'>0$):
\be\label{ozut1}
\bp(x)\;=\;\Re\Big(\bB x+\bB_1x^{2} \Big)+C\,\vec{\gamma}_0|x|^{2}\big(\log|x|^{2}-C_1\big)+\text{O}_2\big(|x|^{\kappa+2-\epsilon'}\big)\:,
\ee
for some constant vectors $\vec{B}_1\in\C^m$ and some real-valued constants $C$, $C_1$. Naturally, the constant $\vec{B}$ remains as in (\ref{locexphi}). 

\subsubsection{Inverted embedded minimal surfaces in Schwarzschild spaces and the proof of Corollary \ref{intromin}}

In this section, we will supposed that our immersion $\bp$ is obtained from inverting an embedded minimal surface. Clearly, $\bp$ is Willmore (away from the singularity at the origin of the unit disk) and it is conformal with respect to the ambient asymptotically Schwarzschild metric of the type (\ref{schwatooo}). Thus all which has been established in section \ref{allld} remains valid. From (\ref{ozut1}), we know that $\bp$ satisfies locally around $z=0$:
\bes
\bp(z,\bar{z})\;=\;\Re\Big(\bB z+\bB_1z^{2}\Big)+C\,\vec{\gamma}_0|z|^{2}\big(\log|z|^{2}-C_1\big)+\text{O}_2\big(|z|^{\kappa+2-\epsilon'}\big)\:,\qquad\forall\:\:\epsilon'>0\:.
\ees
Thus te original immersion $\vec{\xi}=|\bp|^{-2}\bp$ satisfies in particular
\be\label{zut7}
\vec{\xi}(z,\bar{z})\;=\;b^{-2}\Re\Big(\bB \bar{z}^{-1}+\bB_1\bar{z}^{-1}z\Big)+b^{-2}C\,\vec{\gamma}_0\big(\log|z|^{2}-C_1\big)+\text{o}_2(|z|^{\kappa/2})\:,
\ee
where we have used (\ref{stuffonb}) to find that $|\bp|^2\simeq|\bp|^2_g\simeq b^2|z|^2$, with $b^2:=|\Re(\bB)|^2_g+|\Im(\bB)|^2_g$.\\
We are assuming that $\bH^{h}_{\vec{\xi}}=\vec{0}$ away from $z=0$, where, as in Section \ref{reform}, $h$ denotes the ambient metric on $\R^m$ prior to the inversion. Using a formula akin to (\ref{laseulefaute}) with $h$ in place of $g$ shows that
$\Delta\bxi=\text{O}(1)$. Note that $\Delta\Re(z/\bar{z})\simeq|z|^{-2}\gg\text{o}(|z|^{\kappa/2-2})$. From this it follows that $\bB_1$ in the expansion (\ref{zut7}) must be $\vec{0}$. This yields a local expansion of the form
\be\label{zut8}
\vec{\xi}(z,\bar{z})\;=\;\Re\big(\vec{a} \bar{z}^{-1}\big)+\vec{a}_2+b^{-2}C\,\vec{\gamma}_0\log|z|^{2}+\text{O}_2(|z|^{\kappa-\epsilon'})\:,\qquad\forall\:\:\epsilon'>0\:,
\ee
with $\vec{a}:=b^{-2}\Re\bB$.  It is easy to verify that $\vec{a}$ inherits from (\ref{stuffonb}) the properties:
\bes
|\vec{a}_R|_h\,=\,|\vec{a}_I|_h\:\:,\quad
\langle\vec{a}_R,\vec{a}_I\rangle_{h}\,=\,0\:\:,\quad\text{and}\quad \pi_{\bn_h(0)}\vec{a}\,=\,\vec{0}\:.
\ees
Note that we have again used that the metric $h$ defined in (\ref{condAF}) is equivalent to the metric $g$. \\
Because near the origin, $\vec{a}$ is a tangent vector, while, owing to (\ref{zut2}), $\vec{\gamma}_0$ is normal vector, it is not difficult to see that (\ref{zut8}) can be recast as a graph over $\R^2\setminus D_R(0)$:
\bes
(r,\varphi)\,\longmapsto\,\big(r\cos\varphi\,,r\sin\varphi\,,\vec{a}_0+\vec{c}_0\log r+\text{O}_2(r^{\kappa-\epsilon'})\big)\:,\qquad\forall\:\:\epsilon'>0\:,
\ees
in the range $\varphi\in[0,2\pi)$ and $r>R$, for some $R$ chosen large enough, and for some $\R^m$-valued constant vectors $\vec{a}_0$ and $\vec{c}_0$.\\

\bigskip

\renewcommand{\theequation}{A.\arabic{equation}}
\renewcommand{\theTh}{A.\arabic{Th}}
\renewcommand{\theProp}{A.\arabic{Prop}}
\renewcommand{\theLma}{A.\arabic{Lma}}
\renewcommand{\theCo}{A.\arabic{Co}}
\renewcommand{\theRm}{A.\arabic{Rm}}
\renewcommand{\theequation}{A.\arabic{equation}}
\setcounter{equation}{0} 
\reset
\appendix
\section{Appendix}
\subsection{Notational Conventions}

We append an arrow to all the elements belonging to $\R^m$. To simplify the notation, by $\bp\in X(\di)$ is meant $\bp\in X(\di,\R^m)$ whenever $X$ is a function space. Similarly, we write $\nabla\bp\in X(\di)$ for $\nabla\bP\in \mathbb{R}^2\otimes X(\di,\R^m)$.\\[1.5ex]
We let differential operators act on elements of $\R^m$ componentwise. Thus, for example, $\nabla\bp$ is the element of $\R^2\otimes\R^m$ with $\R^m$-valued components $(\px\bp,\py\bp)$. If $S$ is a scalar and $\bR$ an element of $\R^m$, then we let
\begin{eqnarray*}
\bR\cdot\nabla\bP&:=&\big(\bR\cdot\px\bP\,,\,\bR\cdot\py\bP\big)\:\\[1ex]
\nabla^\perp S\cdot\nabla\bP&:=&\px S\,\py\bp\,-\,\py S\,\px\bp\:\\[1ex]
\nabla^\perp\bR\cdot\nabla\bP&:=&\px\bR\cdot\py\bp\,-\,\py\bR\cdot\px\bp\:\\[1ex]
\nabla^\perp\bR\wedge\nabla\bP&:=&\px\bR\wedge\py\bp\,-\,\py\bR\wedge\px\bp\:.
\end{eqnarray*}
Analogous quantities are defined according to the same logic. \\

Two operations between multivectors are useful. The {\it interior multiplication} $\res$ maps a pair comprising a $q$-vector $\gamma$ and a $p$-vector $\beta$ to a $(q-p)$-vector. It is defined via
\bes
\langle \gamma\res\beta\,,\alpha\rangle\;=\;\langle \gamma\,,\beta\wedge\alpha\rangle\:\qquad\text{for each $(q-p)$-vector $\al$.}
\ees
Let $\al$ be a $k$-vector. The {\it first-order contraction} operation $\bul$ is defined inductively through 
\bes
\al\bul\beta\;=\;\al\res\beta\:\:\qquad\text{when $\beta$ is a 1-vector}\:,
\ees
and
\bes
\al\bul(\beta\wedge\gamma)\;=\;(\al\bul\beta)\wedge\gamma\,+\,(-1)^{pq}\,(\al\bul\gamma)\wedge\beta\:,
\ees
when $\beta$ and $\gamma$ are respectively a $p$-vector and a $q$-vector.

\subsection{Miscellaneous Facts}

\subsubsection{Willmore system}

We establish in this section a few general identities. We let $\bp$ be a (smooth) conformal immersion of the unit-disk into $(\R^m,g)$ with $\bp(0)=\vec{0}$. We suppose the metric $g$ satisfies
\be\label{methyppp}
g_{\al\beta}(y)\;=\;\delta_{\al\beta}+\text{O}_2\big(|y|^{\tau}\big)\quad,\quad|y|\ll1\:,
\ee
for some $\tau>0$. As $\bp$ is conformal, the induced metric satisfies
\bes
\tilde{g}_{ij}\;:=\;\big\langle\partial_{x^i}\bp\,,\partial_{x^j}\bp\big\rangle_g\;=\;\text{e}^{2\la}\delta_{ij}\:.
\ees
We will also need the metric $\tilde{g}_0$ induced by pulling back via $\bp$ the Euclidean metric of $\R^m$ on the unit-disk. According to (\ref{methyppp}), one checks that
\bes
(\tilde{g}_0)_{ij}\;=\;\text{e}^{2\la}\big(\delta_{ij}+\text{O}_2(|\bp|^{\tau})\big)\qquad\text{and}\qquad |\tilde{g}_0|\;=\;\text{e}^{4\la}\big(1+\text{O}_2(|\bp|^{\tau})\big)\:.
\ees

For notational convenience, we set $\bej:=\text{e}^{-\la}\pj\bp$. Since $\bp$ is conformal, $\{\bex,\bey\}$ forms an orthonormal basis of the tangent space for the metric $g$. Let $\bn_g:=\star_g(\bex\wedge\bey)$. If $\bV$ is a 1-vector, we find
\begin{eqnarray*}
(\star_g\bn_g)\cdot(\bV\wedge\pj\bp)&=&\text{e}^{\la}(\bex\wedge\bey)\cdot(\bV\wedge\bej)\;=\;\text{e}^{-\la}\big[(\bAe_{1}\cdot\bV )(\tilde{g}_0)_{2j}-(\bAe_{2}\cdot\bV )(\tilde{g}_0)_{1j} \big]\\[1ex]
&=&\text{e}^{\la}\big[(\bAe_{1}\cdot\bV )\delta_{2j}-(\bAe_{2}\cdot\bV )\delta_{1j} \big]+\text{O}_2\big(\text{e}^{\la}|\bp|^{\tau}|\bV|  \big)\\[1ex]
&=&-\,\bV\cdot\partial_{x^{j'}}\bp\,+\text{e}^{\la}|\bV|\text{O}_2(\bp|^{\tau})\:,
\end{eqnarray*}
where
\bes
(\partial_{x^{1'}},\partial_{x^{2'}})\,:=\,(\partial_{x^{2}},-\partial_{x^{1}})\:.
\ees
Hence,
\be\label{idd1}
\left\{\begin{array}{lcr}
(\star_g\bn_g)\cdot(\bV\wedge\nabla\bp)&=&\bV\cdot\nabla^\perp\bp\,+\text{e}^{\la}|\bV|\text{O}_2(\bp|^{\tau})\\[1ex]
(\star_g\bn_g)\cdot(\bV\wedge\nabla^\perp\bp)&=&-\,\bV\cdot\nabla\bp\,+\text{e}^{\la}|\bV|\text{O}_2(\bp|^{\tau})\:.
\end{array}\right.
\ee

\medskip
We choose next an orthonormal basis $\{\bn_\al\}^{m-2}_{\al=1}$ of the normal space such that $\,\{\bex,\bey,\bn_1,\ldots,\bn_{m-2}\}$ is a positive oriented orthonormal basis of $(\R^m,g)$. \\
Recalling the definition of the interior multiplication operator $\res$ given in Section A.1 (understood here for the Euclidean metric in $\R^m$) , it is not hard to obtain
\begin{eqnarray}
(\star_g\bn_g)\bul(\bej\wedge\bek)&=&\text{e}^{-2\la}\Big[(\tilde{g}_0)_{2k}\,\bAe_1\wedge\bAe_{j}-(\tilde{g}_0)_{2j}\,\bAe_1\wedge\bAe_{k}-(\tilde{g}_0)_{1k}\,\bAe_2\wedge\bAe_{j}+(\tilde{g}_0)_{1j}\,\bAe_2\wedge\bAe_{k}  \Big]\nonumber\\[1ex]
&=&\text{O}_2\big(|\bp|^{\tau}  \big)\:,
\end{eqnarray}
and
\begin{eqnarray}
(\star_g\bn_g)\bul(\bn_\al\wedge\bej)&=&\text{e}^{-2\la}\Big[(\tilde{g}_0)_{1j}\,\bn_\al\wedge\bAe_{2}-(\tilde{g}_0)_{2j}\,\bn_\al\wedge\bAe_{1}\Big]+(\bn_\al\cdot\bAe_{2})\,\bAe_1\wedge\bAe_j-(\bn_\al\cdot\bAe_{1})\,\bAe_2\wedge\bAe_j\nonumber\\[1ex]
&=&\delta_{1j}\,\bn_\al\wedge\bAe_{2}-\delta_{2j}\,\bn_\al\wedge\bAe_{1}\,+\text{e}^{\la}|\bV|\text{O}_2(\bp|^{\tau})\:.
\end{eqnarray}
From this one easily deduces for every 1-vector $\bV$, one has
\be\label{idd20}
\left\{\begin{array}{lcl}
(\star_g\bn_g)\bul\big(\bV\wedge\nabla\bp\big)&=&\pi_{\bn_g}\bV\wedge\nabla^\perp\bp\,+\text{e}^{\la}|\bV|\text{O}_2(\bp|^{\tau})\\[1.75ex]
(\star_g\bn_g)\bul\big(\bV\wedge\nabla^\perp\bp\big)&=&-\,\pi_{\bn_g}\bV\wedge\nabla\bp\,+\text{e}^{\la}|\bV|\text{O}_2(\bp|^{\tau})\:.
\end{array}\right.
\ee

There holds furthermore
\begin{eqnarray}
(\bV\wedge\bej)\bul\bei&=&(\bei\cdot\bV)\,\bej\,-\,(\tilde{g}_0)_{ij}\bV\:.
\end{eqnarray}
Hence:
\be\label{idd3}
\left\{\begin{array}{lll}
\big(\bV\wedge\nabla^\perp\bp\big)\bul\nabla^\perp\bp&=&\text{e}^{2\la}\big(\pi_{T_0}\bV-2\bV\big)\,+\text{e}^{2\la}|\bV|\text{O}_2(\bp|^{\tau})\\[1.5ex]
\big(\bV\wedge\nabla\bp\big)\bul\nabla^\perp\bp&\equiv&-\,\big(\bV\cdot\nabla\bp\big)\cdot\nabla^\perp\bp\:.
\end{array}\right.
\ee
As usual, $\pi_{T_0}\bV$ denotes the tangential part of the vector $\bV$ with respect to the Euclidean metric on $\R^m$. 
\bigskip
We are now sufficiently geared to prove

\begin{Lma}\label{identities}
Let $\bp$ be a smooth conformal immersion of the unit-disk into $(\R^m,g)$, with $g$ as above, and let $\bL$ and $\bU$ be two 1-vectors. Suppose that $\pi_{\bn_g}\bU=\bU$ (i.e. $\bU$ is a normal vector). We define $A\in \R^2\otimes\bigwedge^0(\R^m)$ and $\bB\in\R^2\otimes\bigwedge^2(\R^m)$ via
\bes
\left\{\begin{array}{rclll}
A&=&\bL\cdot\nabla\bP\,-\,\bU\cdot\nabla^\perp\bp&&\\[1.5ex]
\bB&=&\bL\wedge\nabla\bP\,-\,\bU\wedge\nabla^\perp\bP\:.
\end{array}\right.
\ees
Then the following identities hold:
\be\label{hyperid}
\left\{\begin{array}{rclll}
A&=&-\,(\star_g\bn_g)\cdot\bB^\perp+\text{e}^{\la}\big(|\bL|+|\bU|\big)\text{O}_2(\bp|^{\tau})&&\\[1.5ex]
\bB&=&-\,(\star_g\bn_g)\bul\bB^\perp\,+\,(\star_g\bn_g)\,A^\perp+\text{e}^{\la}\big(|\bL|+|\bU|\big)\text{O}_2(\bp|^{\tau})\:,
\end{array}\right.
\ee
where $\,\star_g\bn_g:=(\px\bp\wedge\py\bp)/|\px\bp\wedge\py\bp|_g\,$.\\[1ex]
Moreover, we have
\be\label{hyperdel}
A\cdot\nabla^\perp\bp\,+\,\bB\bul\nabla^\perp\bp\;=\;2\text{e}^{2\la}\bU\,+\text{e}^{\la}|\bU|\text{O}_2(\bp|^{\tau})\:.
\ee
\end{Lma}
$\textbf{Proof.}$
The identities (\ref{idd1}) give immediately the required
\bes
(\star_g\bn_g)\cdot\bB^\perp\;=\;-\,\bL\cdot\nabla\bp\,+\,\bU\cdot\nabla^\perp\bp\;=\;-\,A+\text{O}_2\big(\text{e}^{\la}|\bp|^{1-\frac{1}{\theta_0}+\tau}(|\bL|+|\bU|) \big)\:.
\ees
Analogously, using the fact that $\bU$ is a normal vector, the identities (\ref{idd20}) give
\begin{eqnarray*}
(\star_g\bn_g)\bul\bB^\perp&=&-\,\pi_{\bn_g}\bL\wedge\nabla\bp\,+\,\pi_{\bn_g}\bU\wedge\nabla^\perp\bp+\text{e}^{\la}\big(|\bL|+|\bU|\big)\text{O}_2(\bp|^{\tau})\nonumber\\[1ex]
&=&-\,\bB\,+\pi_{T_g}\bL\wedge\nabla\bp+\text{e}^{\la}\big(|\bL|+|\bU|\big)\text{O}_2(\bp|^{\tau})\nonumber\\[1ex]
&=&-\,\bB\,+\,\Big[\big\langle\bL\,,\nabla^\perp\bp\big\rangle_g+\big\langle\bU\,,\nabla\bp\big\rangle_g\Big]\,(\star_g\bn_g)+\text{e}^{\la}\big(|\bL|+|\bU|\big)\text{O}_2(\bp|^{\tau})\nonumber\\[1ex]
&=&-\,\bB\,+\,(\star_g\bn_g)\,A^\perp+\text{e}^{\la}\big(|\bL|+|\bU|\big)\text{O}_2(\bp|^{\tau})\:,
\end{eqnarray*}
which is the second equality in (\ref{hyperid}). \\
In order to prove (\ref{hyperdel}), we will use (\ref{idd3}). Namely,
\begin{eqnarray*}
\bB\,\bul\nabla^\perp\bp&=&-\,\big(\bL\cdot\nabla\bp\big)\cdot\nabla^\perp\bp\,+\,\text{e}^{2\la}\big(\pi_{T_0}\bU-2\bU\big)\,+\text{e}^{2\la}|\bU|\text{O}_2(\bp|^{\tau})\\[1.5ex]
&=&-\,A\cdot\nabla^\perp\bp\,-\,(\bU\cdot\nabla^\perp\bp)\cdot\nabla^\perp\bp\,+\,\text{e}^{2\la}\big(\pi_{T_0}\bU-2\bU\big)\,+\text{e}^{2\la}|\bU|\text{O}_2(\bp|^{\tau})\\[1.5ex]
&=&-\,A\cdot\nabla^\perp\bp\,-\,2\,\text{e}^{2\la}\bU\,+\text{e}^{2\la}|\bU|\text{O}_2(\bp|^{\tau})\:,
\end{eqnarray*}
which gives the desired identity.\\[-1.5ex]

$\hfill\blacksquare$ \\

\subsection{Nonlinear and weighted elliptic results}

In \cite{BR} and in \cite{Ber1}, analogous versions of the following three results are proved. The versions stated here are slightly different than those appearing in the aforementioned references. Only very minor modifications are needed; details are left to the reader.

\begin{Prop}\label{appdev}
Let $u\in C^2(\di\setminus\{0\})$ and $V\in C^1(\di\setminus\{0\})$ satisfy the equation
\bes\label{vdb}
\text{div}\,\big(\nabla u(x)+V(x,u)\big)\;=\;0\qquad\text{in}\:\: \di\setminus\{0\}\:.
\ees
Assume that for some integer $a\ge1$, and some $p\in(1,\infty)$ there holds
\bes
|x|^aV\:\:,\:\:|x|^{a-1}u\:\:\in\, L^p(\di)\:.
\ees
Then we have
\bes
|x|^{a}\nabla u\,\in\,L^p(\di)\:.
\ees
\end{Prop}

\medskip

\begin{Prop}\label{adams}
Let $u\in W^{1,2}(\di)\cap C^1(\di\setminus\{0\})$ satisfy the equation
\bes
-\,\Delta u\:=\:\nabla b\cdot\nabla^\perp u+\text{div}\,(w)\qquad\quad\text{in}\:\:\:\:\di\:,
\ees
where $w\in L^{2+\eta}(\di)$, for some $\eta>0$. Moreover, suppose
\bes
\Vert\nabla b\Vert_{L^2(\di)}\;<\;\eps_0\:,
\ees
for some $\eps_0$ chosen to be ``small enough". Then
\bes
\nabla u\,\in\,L^{2+\eta}(\di)\:.
\ees
\end{Prop}

\medskip

\begin{Prop}\label{CZ}
Let $u\in C^2(\di\setminus\{0\})$ solve
\bes
\Delta u(x)\;=\;\mu(x)f(x)\qquad\text{in\:\:\:} \di\:,
\ees
where $f\in L^p(\di)$ for some $p>2$. The weight $\mu$ satisfies
\bes
|\mu(x)|\;\simeq\;|x|^{b}\qquad\quad\text{for some}\:\:\: b\in\mathbb{N}\:.
\ees
Then 
\begin{itemize}
\item[(i)] there holds\footnote{$\overline{x}$ is the complex conjugate of $x$. We parametrize $\di$ by $x=x_1+i\,x_2$, and then $\overline{x}:=x_1-i\,x_2$. In this notation, $\nabla u$ in (\ref{stmt0}) is understood as $\partial_{x_1}u+i\,\partial_{x_2}u$.}
\be\label{stmt0}
\nabla u(x)\;=\;P(\overline{x})\,+\,|\mu(x)|\,T(x)\:,
\ee
where $P(\overline{x})$ is a complex-valued polynomial of degree at most $b$, and near the origin $\,T(x)=\text{O}\big(|x|^{1-\frac{2}{p}-\epsilon'}\big)$ for every $\epsilon'>0$. \\[-1.5ex]
\item[(ii)] furthermore, if $\,\mu\in C^1(\di\setminus\{0\})$ and if
\bes\label{hypw2}
|x|^{1-b}\,\nabla\mu(x)\,\in\,\bigcap_{s<\infty}L^{s}(\di)\:,
\ees
there holds
\bes\label{stmt2}
\nabla^2 u(x)\;=\;\nabla P(\overline{x})\,+\,|\mu(x)|\,Z(x)\:,
\ees
where $P$ is as in (i), and 
\bes
Z\,\in\,L^{p-\epsilon'}(\di,\C^2)\qquad\quad\forall\:\:\epsilon'>0\:.
\ees
\end{itemize}
\end{Prop}

\end{document}